\documentclass[12pt]{amsart}
\usepackage{amsmath,amssymb,latexsym,cancel,rotating}

\numberwithin{equation}{section}

\newtheorem{theorem}{Theorem}[section]
\newtheorem{proposition}[theorem]{Proposition}
\newtheorem{conjecture}[theorem]{Conjecture}

\newtheorem{lemma}[theorem]{Lemma}

\theoremstyle{definition}

\newtheorem{remark}[theorem]{Remark}
\newtheorem{example}[theorem]{Example}
\newtheorem{definition}[theorem]{Definition}

\input xy
\xyoption{poly}
\xyoption{2cell}
\xyoption{all}

\def\A{\mathcal{A}}

\def\NN{\Bbb N}

\def\ZZ{\mathbb Z}

\def\PP{\mathbb{P}}

\def\ll{\mathbb{L}}

\def\ZZ{\mathbb{Z}}

\renewcommand{\eqref}[1]{{\rm (\ref{#1})}}

\begin{document}

\title[Generalized quantum cluster algebras]
{Generalized quantum cluster algebras:  the Laurent phenomenon and upper bounds}

\author{Liqian Bai, Xueqing Chen, Ming Ding and Fan Xu}
\address{Department of Applied Mathematics, Northwestern Polytechnical University, Xi'an, Shaanxi 710072, P.R. China}
\email{bailiqian@nwpu.edu.cn (L.Bai)}
\address{Department of Mathematics,
 University of Wisconsin-Whitewater\\
800 W. Main Street, Whitewater, WI. 53190. USA}
\email{chenx@uww.edu (X.Chen)}
\address{School of Mathematics and Information Science\\
Guangzhou University, Guangzhou 510006, P.R.China}
\email{m-ding04@mails.tsinghua.edu.cn (M.Ding)}
\address{Department of Mathematical Sciences\\
Tsinghua University\\
Beijing 100084, P.~R.~China} \email{fanxu@mail.tsinghua.edu.cn(F.Xu)}

\thanks{Liqian Bai was supported by NSF of China (No. 11801445) and the Natural Science Foundation of Shaanxi Province (NO. 2020JQ-116), Ming Ding was supported by NSF of China (No. 11771217) and Fan Xu was supported by NSF of China (No. 12031007).}




\maketitle

\begin{abstract}
Generalized quantum cluster algebras  introduced in  \cite{BCDX} are  quantum deformation of  generalized cluster algebras of geometric types. In this paper, we  prove that the Laurent phenomenon holds in these generalized quantum cluster algebras. We also show that  upper bounds coincide with the corresponding generalized quantum upper cluster algebras under the ``coprimality" condition.
\end{abstract}


\section{Introduction}
Cluster algebras were invented by Fomin and Zelevinsky \cite{ca1,ca2} in order
to set up an algebraic framework for studying the total positivity in algebraic groups and the dual canonical bases in quantized enveloping algebras. The Laurent phenomenon holds for all cluster variables in  cluster algebras \cite{ca1}, i.e., all cluster variables belong to an intersection of certain (may be infinitely many) rings of Laurent polynomials  which is called  an upper cluster algebra by Berenstein, Fomin and Zelevinsky in \cite{bfz}. They also introduced the upper bound as a certain finite intersection of rings of Laurent polynomials.
These three algebras satisfy the following relations:
\begin{equation*}
  \text{cluster algebra}\subseteq \text{upper cluster algebra}\subseteq \text{upper bound}.
\end{equation*}
They proved that  for a cluster algebra possessing an acyclic and coprime seed,  all three algebras coincide.

As a natural generalization,
Chekhov and Shapiro \cite{CS} introduced the concept of  generalized cluster algebras  in which the  binomial exchange
relations for cluster variables of  cluster algebras are replaced by the polynomial exchange relations for those cluster variables of generalized cluster algebras. Generalized cluster algebras are studied not only in a similar way as cluster algebras~\cite{nak2, CL1, CL2} but also on the scattering diagrams setting~\cite{Mou}.
Gekhtman, Shapiro and Vainshtein \cite{gsv18} proved that generalized upper cluster algebras
coincide with upper bounds under certain coprimality
conditions.  In ~\cite{BCDX-1}, we proved that generalized cluster algebras coincide with upper bounds under the conditions of
acyclicity and coprimality.  Du and Li  \cite{DL} proved that the above results also hold in Laurent phenomenon algebras under some certain conditions.

The quantum counterparts of cluster algebras called quantum cluster algebras were introduced by Berenstein and Zelevinsky \cite{BZ2005} in order to study the dual canonical bases. As we known, all above mentioned results on cluster algebras have been extended to the quantum setting such as the Laurent phenomenon and the coincidence of quantum upper cluster algebras and upper bounds.

Generalized quantum cluster algebras introduced in \cite{BCDX} naturally extend the definition of quantum cluster algebras and  can  be considered as quantum analogues of  generalized cluster algebras of geometric types \cite{CS}. We proved in~\cite{BCDX} that the Laurent phenomenon holds only in these generalized quantum cluster algebras of rank two. In this paper, we obtain the Laurent phenomenon in  all generalized quantum cluster algebras of geometric types, i.e., all generalized cluster variables belong to an intersection of certain (infinitely many) rings of Laurent polynomials which is called generalized quantum upper cluster algebras. We study the upper bound which is the intersection of rings of Laurent polynomials in a finite collection of generalized clusters, through mutations of a fixed cluster along all directions.
Upper bounds have their own interests in cluster theory \cite{bfz,gsv18}.
In order to prove that the upper bound is independent of the choice of the initial cluster, we need the ``coprimality" condition.
In cluster algebras and quantum cluster algebras, the ``coprimality" condition holds if the extended exchange matrix $\widetilde{B}$ has full rank.
Unfortunately, in generalized (quantum) cluster algebras, the ``coprimality" condition cannot be deduced from the full rank of $\widetilde{B}$ (see \cite[Lemma 4.3]{gsv18}).
 Then we show that upper bounds coincide with the corresponding generalized quantum upper cluster algebras under the requirement of the ``coprimality" condition. As a byproduct, the Laurent phenomenon can be implied directly in these generalized quantum cluster algebras. 

The paper is organized as follows. In Section 2, we recall the definition of the generalized quantum cluster algebras, and introduce the notion of the generalized quantum upper cluster algebras.
Section 3 follows the arguments in \cite{gy} with necessary modifications to provide the first main result of this paper concerning the Laurent phenomenon of the generalized quantum cluster algebras. Lastly, in Section 4, under the ``coprimality" condition, we prove that  upper bounds are invariant under mutations and thus coincide with the corresponding generalized quantum upper cluster algebras by using the methods developed in \cite{bfz,BZ2005}.

\section{Generalized quantum cluster algebras}
For simplicity, let $[i,j]:=\{i,i+1,\ldots,j\}$ for $i,j\in\ZZ$ and $i<j$. At first, we recall the definition of the generalized quantum cluster algebras in \cite{BCDX}.

\begin{definition}
Let $m$ and $n$ be positive integers and $m\geq n$. Let $\widetilde{B}=(b_{ij})$ be an $m\times n$ integer matrix and $\Lambda=(\lambda_{ij})$ a skew-symmetric $m\times m$ integer matrix.
The pair $(\Lambda,\widetilde{B})$ is called compatible if there exists a diagonal matrix $D=\text{diag}(\widetilde{d}_1,\widetilde{d}_2,\ldots,\widetilde{d}_n)$, where $\widetilde{d}_i\in\ZZ_{>0}$ for all $i\in[1,n]$, such that
\begin{equation*}
\Lambda\widetilde{B}=-
  \begin{bmatrix}
  D \\
  0
  \end{bmatrix}.
\end{equation*}
\end{definition}

Note that if $(\Lambda,\widetilde{B})$ is a compatible pair, then $\widetilde{B}=
  \begin{bmatrix}
  B \\
  C
  \end{bmatrix}$
has full rank and the matrix $DB$ is skew-symmetric (\cite[Theorem 3.3]{BZ2005}). The matrix $B$ is called the principal part of $\widetilde{B}$. The skew-symmetric matrix $\Lambda$ also gives the skew-symmetric bilinear form $\Lambda:\ZZ^{m}\times\ZZ^{m}\rightarrow\ZZ$ defined by $\Lambda(a,b)=a^T\Lambda b$ for any column vectors $a,b\in\ZZ^{m}$.

For each $k\in[1,n]$, let $d_k$ be a positive integer such that $\frac{b_{lk}}{d_k}$ are integers for all $l\in[1,m]$. The $k$-th column of $\widetilde{B}$ is denoted by $b^{k}$. Let $\beta^{k}:=\frac{1}{d_k}b^{k}$.

If the pair $(\Lambda,\widetilde{B})$ is compatible, using the fact that $\sum\limits_{k=1}^{m}b_{kj}\lambda_{ki}=\delta_{ij}\widetilde{d}_{j}$, it follows that for the standard column vector $e_l$ we have
\begin{equation}\label{equ1}
\Lambda(\beta^{k},e_l)=\left\{
  \begin{aligned}
    \widetilde{d}_kd_{k}^{-1}&,~\text{if}~l=k;  \\
    0{\hskip 0.5cm}&,~\text{otherwise}.
  \end{aligned}
  \right.
\end{equation}

Let $(\Lambda,\widetilde{B})$ be a compatible pair and the sign $\varepsilon\in\{\pm1\}$. Recall that the function $[x]_+=x$ if $x\geq0$, and $[x]_+=0$ if $x\leq0$. For each $i\in[1,n]$, the matrices $\widetilde{B}^\prime=\mu_i(\widetilde{B})$ and $\Lambda^\prime=\mu_i(\Lambda)$ are defined by
$$
\widetilde{B}^\prime= E_{\varepsilon}^{\widetilde{B}}\widetilde{B}F_{\varepsilon}^{\widetilde{B}}~\text{and}~ \Lambda^\prime= (E_{\varepsilon}^{\widetilde{B}})^T\Lambda E_{\varepsilon}^{\widetilde{B}},
$$
where the $m\times m$ matrix $E_{\varepsilon}^{\widetilde{B}}=E_{\varepsilon}$ is defined as follows
\begin{align*}
(E_{\varepsilon})_{kl}=\left\{
  \begin{aligned}
  \delta_{kl}{\hskip 0.7cm}&,~\text{if}~l\neq i;\\
  -1{\hskip 0.7cm}&,~\text{if}~k=l=i;\\
  [-\varepsilon b_{ki}]_+&,~\text{if}~k\neq l=i,
  \end{aligned}
\right.
\end{align*}
and the $n\times n$ matrix $F_{\varepsilon}^{\widetilde{B}}=F_{\varepsilon}$ is defined by
\begin{align*}
(F_{\varepsilon})_{kl}=\left\{
  \begin{aligned}
  \delta_{kl}{\hskip 0.3cm}&,~\text{if}~k\neq i;\\
  -1{\hskip 0.3cm}&,~\text{if}~k=l=i;\\
  [\varepsilon b_{il}]_+&,~\text{if}~k=i\neq l.
  \end{aligned}
\right.
\end{align*}

By \cite[Proposition 3.4]{BZ2005}, it follows that the pair $(\Lambda^{\prime},{\widetilde{B}}^{\prime})$ is also compatible and independent of the choice of $\varepsilon$. By \cite[Proposition 3.6]{BZ2005}, we know that $\mu_i$ is an involution. The pair $(\Lambda^{\prime},{\widetilde{B}}^{\prime})$ is said to be the mutation of $(\Lambda, {\widetilde{B}})$ in direction $i$.

Let $\mathbb{Z}[q^{\pm\frac{1}{2}}]\subset \mathbb{Q}(q^{\frac{1}{2}})$ be the ring of integer Laurent polynomials in a formal variable $q^{\frac{1}{2}}$.

\begin{definition}
The quantum torus $\mathcal{T}=\mathcal{T}(\Lambda)$ is the $\ZZ[q^{\frac{1}{2}}]$-algebra with the distinguished basis $\{X(c)~|~c\in\ZZ^{m}\}$ and the multiplication defined as
\begin{equation}\label{qcommutation}
X(c)X(d)=q^{\frac{1}{2}\Lambda(c,d)}X(c+d),
\end{equation}
where $c,d\in\ZZ^{m}$. Let $\mathcal{F}$ denote the skew-field of fractions of $\mathcal{T}$.
\end{definition}

It follows that $X(c)X(d)=q^{\Lambda(c,d)}X(c)X(d)$, $X(0)=1$ and $X(-c)=X(c)^{-1}$ for $c,d\in\ZZ^{m}$. Let $X_k:=X(e_k)$ for $k\in[1,m]$, we have that
$$
X(c)=q^{\frac{1}{2}\sum\limits_{l< k}c_lc_k\lambda_{kl}}X_{1}^{c_1}X_{2}^{c_2}\ldots X_{m}^{c_m}
$$
for $c=(c_1,c_2,\ldots,c_m)^{T}\in\ZZ^{m}$.

For $k\in[1,n]$, let
$$
\mathbf{h}_k:=\{h_{k,0}(q^{\frac{1}{2}}),h_{k,1}(q^{\frac{1}{2}}),\ldots,h_{k,d_k}(q^{\frac{1}{2}})\},
$$
where $h_{k,l}(q^{\frac{1}{2}})\in\ZZ[q^{\pm\frac{1}{2}}]$ are Laurent polynomials satisfying that $h_{k,l}(q^{\frac{1}{2}})= h_{k,d_k-l}(q^{\frac{1}{2}})$ and $h_{k,0}(q^{\frac{1}{2}})= h_{k,d_k}(q^{\frac{1}{2}})=1$. Let $\mathbf{h}:=(\mathbf{h}_1,\mathbf{h}_2,\ldots,\mathbf{h}_n)$.

\begin{definition}
Let $(\Lambda,\widetilde{B})$ be a compatible pair. The quadruple $(X,\mathbf{h},\Lambda,\widetilde{B})$ is called a quantum seed. For $i\in[1,n]$, the new quadruple $$(X^\prime,\mathbf{h}^\prime,\Lambda^\prime,\widetilde{B}^\prime):=\mu_i(X,\mathbf{h},\Lambda,\widetilde{B})$$ is defined as follows
\begin{equation}\label{clustervarialbemutaion}
X^\prime(e_k)=\mu_i(X(e_k))=\left\{
  \begin{aligned}
    X(e_k), {\hskip 6.5cm}&~\text{if}~k\neq i;  \\
    \sum\limits_{r=0}^{d_i}h_{i,r}(q^{\frac{1}{2}})X(r[\beta^i]_+ +(d_i-r)[-\beta^{i}]_+ -e_i), &~\text{if}~k=i,
  \end{aligned}
  \right.
\end{equation}
and
$$
\mathbf{h}^\prime=\mu_i(\mathbf{h})=\mathbf{h},~\Lambda^\prime=\mu_i(\Lambda) \text{ and }~\widetilde{B}^\prime=\mu_i(\widetilde{B}).
$$
We say that the quadruple $\mu_i(X,\mathbf{h},\Lambda,\widetilde{B})$ is obtained from $(X,\mathbf{h},\Lambda,\widetilde{B})$ by the mutation in direction $i$. For simplicity, write $X^\prime_k:=X^\prime(e_k)$.
\end{definition}

By \cite[Proposition 3.6 and Proposition 3.7]{BCDX}, the quadruple $(X^\prime,\mathbf{h}^\prime,\Lambda^\prime,\widetilde{B}^\prime)$ is a quantum seed and $\mu_i$ is an involution. The quantum seed $(X^{\prime}, \textbf{h}^{\prime}, \Lambda^{\prime}, \widetilde{B}^{\prime})$ is said to be mutation-equivalent to $(X, \textbf{h}, \Lambda, \widetilde{B})$,
if $(X^{\prime}, \textbf{h}^{\prime}, \Lambda^{\prime}, \widetilde{B}^{\prime})$ can be obtained from
$(X, \textbf{h}, \Lambda, \widetilde{B})$ by a sequence of seed mutations, i.e.,
$$
(X^{\prime}, \textbf{h}^{\prime}, \Lambda^{\prime}, \widetilde{B}^{\prime})=\mu_{i_t}(\ldots(\mu_{i_1}(X, \textbf{h}, \Lambda, \widetilde{B}))\ldots)
$$
for some $1\leq i_1,\ldots,i_t\leq n$, which is denoted by
$$
(X^{\prime}, \textbf{h}^{\prime}, \Lambda^{\prime}, \widetilde{B}^{\prime})\sim (X, \textbf{h}, \Lambda, \widetilde{B}).
$$

The set $\{X^{\prime}_1,\ldots,X^{\prime}_n\}$ (resp. $\{X^{\prime}_1,\ldots,X^{\prime}_m\}$) is called the cluster (resp. an extended cluster) of the quantum seed $(X^{\prime}, \textbf{h}^{\prime}, \Lambda^{\prime}, \widetilde{B}^{\prime})$. The elements $X^{\prime}_k$ are called cluster variables for $k\in[1,n]$ and $X^{\prime}_l$ frozen variables for $l\in[n+1,m]$. In fact $X^{\prime}_l=X_l$ for all $l\in[n+1,m]$.

\begin{definition}\label{def of qgca}
The generalized quantum cluster algebra $\mathcal{A}(X, \textbf{h}, \Lambda, \widetilde{B})$ with the initial quantum seed $(X, \textbf{h}, \Lambda, \widetilde{B})$, is the $\ZZ[q^{\pm\frac{1}{2}}]$-subalgebra of $\mathcal{F}$ generated by the cluster variables $X^{\prime}_1,X^{\prime}_2,\ldots,X^{\prime}_n$ and the frozen variables $X_{n+1}^{\pm1}, X_{n+2}^{\pm1},\ldots,X_{m}^{\pm1}$ from all quantum seeds $(X^{\prime}, \textbf{h}^{\prime}, \Lambda^{\prime}, \widetilde{B}^{\prime})$ which are mutation-equivalent to $(X, \textbf{h}, \Lambda, \widetilde{B})$.
\end{definition}

For convenience, let $\mathbb{P}$ denote the multiplicative group generated by $X_{n+1},\ldots,X_{m}$ and $q^{\frac{1}{2}}$, and $\mathbb{ZP}$  the ring of the Laurent polynomials in $X_{n+1},\ldots, X_m$ with coefficients in $\ZZ[q^{\pm\frac{1}{2}}]$.

\begin{definition}
The generalized quantum upper cluster algebra $\widetilde{\mathcal{U}}(X, \textbf{h}, \Lambda, \widetilde{B})$ is defined as:
\begin{align*}
\widetilde{\mathcal{U}}(X, \textbf{h}, \Lambda, \widetilde{B}) :=
\bigcap\limits_{(X^{\prime}, \textbf{h}^{\prime}, \Lambda^{\prime}, \widetilde{B}^{\prime})\sim (X, \textbf{h}, \Lambda, \widetilde{B}) }\ZZ\mathbb{P}[(X^{\prime}_1)^{\pm1},\ldots,(X^{\prime}_n)^{\pm1}].
\end{align*}
\end{definition}


\section{The Laurent Phenomenon}
In this section, we show the following main result that the Laurent phenomenon also holds in generalized quantum cluster algebras.
\begin{theorem}[Laurent phenomenon]\label{laurentphenomenon}
Let $(X, \textbf{h}, \Lambda, \widetilde{B})$ be a quantum seed. Then we have
$$\A(X, \textbf{h}, \Lambda, \widetilde{B})  \subseteq \widetilde{\mathcal{U}}(X, \textbf{h}, \Lambda, \widetilde{B}).$$
\end{theorem}

To prove Theorem \ref{laurentphenomenon}, we need some lemmas.

Given a quantum seed $(X,\textbf{h},\Lambda_0,\widetilde{B}_0)$. We define
\begin{itemize}
  \item[] $(X^{(0)},\textbf{h},\Lambda_0,\widetilde{B}_0) =(X,\textbf{h},\Lambda_0,\widetilde{B}_0)$,
  \item[] $(X^{(1)},\textbf{h},\Lambda_1,\widetilde{B}_1) =\mu_i( X^{(0)},\textbf{h},\Lambda_0,\widetilde{B}_0 )$,
  \item[] $(X^{(2)},\textbf{h},\Lambda_2,\widetilde{B}_2) =\mu_j(X^{(1)},\textbf{h},\Lambda_1,\widetilde{B}_1)$,
  \item[] $(X^{(3)},\textbf{h},\Lambda_3,\widetilde{B}_3) =\mu_i(X^{(2)},\textbf{h},\Lambda_2,\widetilde{B}_2)$.
\end{itemize}

For simplicity, write $x:=X_i$, $y:=X_j=X^{(1)}_j$, $z:=X^{(1)}_i=X^{(2)}_i$, $u:=X^{(2)}_j=X^{(3)}_j$ and $v:=X^{(3)}_i$.

Let $P:=xz$, $Q:=yu$ and $R:=zv$. Thus $P$ is a polynomial in $X_k$ for $k\in [1,m]\setminus\{i\}$. We can write $P(y):=P$ as a polynomial in $y$ with coefficients in left hand being noncommutative polynomials in $X_k$ for $k\in [1,m]\setminus\{i,j\}$. Similarly, $Q$ (resp. $R$) is a polynomial in $X^{(1)}_k$ (resp. $X^{(2)}_l$) for $k\in [1,m]\setminus\{j\}$ (resp. for $l\in [1,m]\setminus\{i\}$) and $Q(z):=Q$ (resp. $R(u):=R$).

For each $t\in [0,3]$, let
$$
\mathbb{L}_t:=\ZZ\mathbb{P}[ (X^{(t)}_1)^{\pm1},\ldots,(X^{(t)}_{n})^{\pm1}],
$$
which is a Laurent polynomial ring. Similarly, write
$$
\mathbb{L}_{0}^{\vee}:=\ZZ\mathbb{P}[X_1^{\pm1},\ldots,X_{i-1}^{\pm1},X_{i+1}^{\pm1},\ldots,X_{n}^{\pm1}] \text{ and } \mathbb{L}_{0}^{\vee\vee}:=\ZZ\mathbb{P}[X_k^{\pm1}~|~k\in[1,n]\setminus\{i,j\}].$$
For $t\in[0,3]$, the $kl$-th entry of $\widetilde{B}_t$ (resp. $\Lambda_t$) is denoted by $b_{t,kl}$ (resp. $\lambda_{t,kl}$). For simplicity, let $\beta_{t,kl}:=\frac{b_{t,kl}}{d_l}$.

The following result is similar to~\cite[Lemma 2.16]{gy}.
\begin{lemma}\label{coprimelem}
Let $N,N_1\in\ZZ_{>0}$. Let
$$
V_1(X(f_1))=\Big(\sum\limits_{r_1=0}^{s_1}\xi_{1,r_1}X(r_1f_1)\Big)\Big(\sum\limits_{r_2=0}^{s_2}\xi_{2,r_2}X(r_2f_1)\Big) \cdots \Big(\sum\limits_{r_N=0}^{s_N}\xi_{N,r_N}X(r_Nf_1)\Big)
$$
and
$$
V_2(X(f_2))=\Big(\sum\limits_{r^{\prime}_{1}=0}^{t_1}\zeta_{1,r^{\prime}_{1}}X(r^{\prime}_{1}f_2)\Big) \Big(\sum\limits_{r^{\prime}_{2}=0}^{t_2}\zeta_{2,r^{\prime}_{2}}X(r^{\prime}_{2}f_2)\Big)\cdots \Big(\sum\limits_{r^{\prime}_{N_1}=0}^{t_{N_1}}\zeta_{N_1,r^{\prime}_{N_1}}X(r^{\prime}_{N_1}f_2)\Big)
$$
be two elements in $\mathbb{L}_{0}^{\vee\vee}\subseteq\mathcal{T}$, where $f_1,~f_2\in\sum\limits_{k\in[1,m] \setminus\{i,j\}}\ZZ e_k$ satisfying that
\begin{equation*}
\Lambda(f_1,e_k)=\Lambda(f_2,e_k)=0
\end{equation*}
for every $k\in[1,m]\setminus\{i,j\}$, coefficients $\xi_{l_1,r_{l_1}},\zeta_{l_2,r^{\prime}_{l_1}}\in\ZZ[q^{\pm\frac{1}{2}}]$ and $\xi_{l_1,s_{l_1}}=\zeta_{l_2,t_{l_2}}=1$ for $l_1\in[0,N],~l_2\in[0,N_1]$. If $f_1$ and $f_2$ are $\ZZ$-linearly independent, then $V_1(X(f_1))$ and $V_2(X(f_2))$ are coprime in the center of $\mathbb{L}_{0}^{\vee\vee}$.

\end{lemma}

\begin{proof}
First of all, we can assume that $f_1$ and $f_2$ are nonzero. If $\xi_{1,0}=0$, then
$$
\sum\limits_{r_1=0}^{s_1}\xi_{1,r_1}X(r_1f_1)=\Big(\sum\limits_{r_1=1}^{s_1}\xi_{1,r_1}X((r_1-1)f_1)\Big) X(f_1)
$$
and $X(f_1)$ is an invertible element in $\ll_0$. Hence we can replace $\sum\limits_{r_1=0}^{s_1}\xi_{1,r_1}X(r_1f_1)$ with $\sum\limits_{r_1=1}^{s_1}\xi_{1,r_1}X((r_1-1)f_1)$. Therefore we can assume that $\xi_{0,0},\xi_{1,0},\ldots,\xi_{N,0} $ and $\zeta_{0,0},\zeta_{1,0},\ldots, \zeta_{N_1,0}\in\ZZ[q^{\pm\frac{1}{2}}]\setminus\{0\}$.

Let $C(\mathbb{L}_{0}^{\vee\vee})$ denote the center of $\mathbb{L}_{0}^{\vee\vee}$. Note that $C(\mathbb{L}_{0}^{\vee\vee})=\sum\limits_{g\in S}\ZZ[q^{\pm\frac{1}{2}}]X(g)$, where
$$
S=\Big\{g\in\sum\limits_{k\in[1,m] \setminus\{i,j\}}\ZZ e_k~|~\Lambda(g,e_k)=0~\text{for~all~}k\in[1,m]\setminus\{i,j\}\Big\}.
$$
It is clear that $V_1(X(f_1)), V_2(X(f_2)) \in C(\mathbb{L}_{0}^{\vee\vee})$.

Let $F(\ZZ[q^{\pm\frac{1}{2}}])$ denote the fraction field of $\ZZ[q^{\pm\frac{1}{2}}]$ and $\widetilde{F}$ be the algebraic closure of $F(\ZZ[q^{\pm\frac{1}{2}}])$. It follows that $V_1(X(f_1))$ and $V_2(X(f_2))$ are in the center of $\widetilde{F} [X_{1}^{\pm1},\ldots,X_{m}^{\pm1}~|~X_iX_j=q^{\lambda_{ij}}X_jX_i]$.

There exists a positive integer $a_1$ (resp. $a_2$) such that $f_1=a_1f_{1}^{\prime}$ (resp. $f_2=a_2f_{2}^{\prime}$) and the great common divisor of the entries of $f_{1}^{\prime}$ (resp. $f_{2}^{\prime}$) is $1$. The irreducible factors of $V_1(X(f_1))$ and $V_2(X(f_2))$ are of the forms $X(f_{1}^{\prime})+\xi$ and $X(f_{2}^{\prime})+\zeta$, respectively, for some $\xi, \zeta\in\mathbb{L}_{0}^{\vee\vee}$. Using the fact that $f_1$ and $f_2$ are $\ZZ$-linearly independent, we conclude that $f_{1}^{\prime}\neq\pm f_{2}^{\prime}$ and $(X(f_{1}^{\prime})+\xi, X(f_{2}^{\prime})+\zeta)=1$. Hence $V_1(X(f_1))$ and $V_2(X(f_2))$ are coprime elements in the center of $\mathbb{L}_{0}^{\vee\vee}$.
\end{proof}

For $r\in\ZZ$ and $b\in\mathbb{L}_{0}\setminus\{0\}$, the element $cx^{r}\in\mathbb{L}_0$ is said to the leading term of $b$ for some $c\in\mathbb{L}_{0}^{\vee}$ if
$$
b-cx^{r}\in\bigoplus\limits_{r^{\prime}<r}\mathbb{L}_{0}^{\vee}x^{r^{\prime}}.
$$

\begin{lemma}\label{lemL}
The element
$
L:=z^{-1}R(y^{-1}Q(0))zP^{-1}
$
is a monomial in $X_l$ and $y^{-1}$ for $l\in[1,m]\setminus\{i,j\}$.

\end{lemma}
\begin{proof}
When $b_{1,ij}=0$, we have that $Q(0)=Q(z)$. Note that
\begin{align*}
&z^{-1}R(y^{-1}Q(0))=v \\
=&\sum\limits_{r=0}^{d_i}h_{i,r}(q^{\frac{1}{2}})X^{(2)} (r[\beta_{2}^i]_+ +(d_i-r)[-\beta_{2}^i]_+ -e_i)\\
=& \sum\limits_{r=0}^{d_i}h_{i,r}(q^{\frac{1}{2}})X^{(1)} (r[\beta_{2}^i]_+ +(d_i-r)[-\beta_{2}^i]_+ -e_i)
\end{align*}
and
\begin{align*}
Pz^{-1}=x=\sum\limits_{r=0}^{d_i}h_{i,r}(q^{\frac{1}{2}})X^{(1)}(r[\beta_{1}^i]_+ +(d_i-r)[-\beta_{1}^i]_+-e_i).
\end{align*}
We conclude that $L=z^{-1}R(y^{-1}Q(0))zP^{-1}=1$ since $b_{1}^{i}=b_{2}^{i}$.

When $b_{1,ij}>0$, by direct calculations, we obtain that
\begin{align*}
&y^{-1}Q(0)=X^{(1)}([-b_{1}^{j}]_{+}-e_j), \\
&E_{+}^{\widetilde{B}_1}[\beta_{2}^{i}]_+=[\beta_{2}^{i}]_+ -2\beta_{2,ji}e_j +\beta_{2,ji}[-b_{1}^{j}]_+,\\
&E_{+}^{\widetilde{B}_1}[-\beta_{2}^{i}]_+=[-\beta_{2}^{i}]_+, \\
&E_{+}^{\widetilde{B}_1}\beta_{2}^{i}=\beta_{2}^{i} -2\beta_{2,ji}e_j +\beta_{2,ji}[-b_{1}^{j}] =\beta_{1}^{i}, \\
&E_{+}^{\widetilde{B}_1}e_j=[-b_{1}^{j}]_+-e_j.
\end{align*}

It follows that
\begin{align*}
&z^{-1}R(y^{-1}Q(0)) \\
= &X^{(1)}(E_{+}^{\widetilde{B}_1}[b_{2}^{i}]_+)\Big( \sum\limits_{r=0}^{d_i}h_{i,r}(q^{\frac{1}{2}}) q^{-\frac{1}{2}\Lambda_2(r[\beta_{2}^{i}]_{+}+(d_i-r)[-\beta_{2}^{i}]_+,-e_i)}\\
&\cdot q^{-\frac{1}{2}\Lambda_2(r[\beta_{2}^{i}]_{+}+(d_i-r)[-\beta_{2}^{i}]_+,r\beta_{2,ji}e_j)} \cdot q^{\frac{1}{2}\Lambda_1(r[\beta_{2}^{i}]_{+}+(d_i-r)[-\beta_{2}^{i}]_{+}-r\beta_{2,ji}e_j ,r\beta_{2,ji}[-b_{1}^{j}]_{+}-r\beta_{2,ji}e_j)}\\
&\cdot q^{\frac{1}{2}\Lambda_1(-E_{+}^{\widetilde{B}_1}[b_{2}^{i}]_{+},r[\beta_{2}^{i}]_{+}+(d_i-r)[-\beta_{2}^{i}]_{+} -2r\beta_{2,ji}e_j +r\beta_{2,ji}[-b_{1}^{j}]_+)}X^{(1)}((r-d_i)\beta_{1}^{i})\Big)z^{-1}\\
=&X^{(1)}(E_{+}^{\widetilde{B}_1}[b_{2}^{i}]_+)\Big(\sum\limits_{r=0}^{d_i}h_{i,r} (q^{\frac{1}{2}})q^{-\frac{1}{2}\Lambda_2(r[\beta_{2}^{i}]_{+} +(d_i-r)[-\beta_{2}^{i}]_+,-e_i)}X^{(1)}((r-d_i)\beta_{1}^{i}) \Big)z^{-1}\\
=&q^{-\frac{1}{2}\Lambda_2([b_{2}^{i}]_+,-e_i)}X^{(1)}(E_{+}^{\widetilde{B}_1}[b_{2}^{i}]_+) \Big(\sum\limits_{r=0}^{d_i}h_{i,r}(q^{\frac{1}{2}}) q^{\frac{1}{2}(r-d_i)\widetilde{d}_id_{i}^{-1}}X^{(1)}((r-d_i)\beta_{1}^{i}) \Big)z^{-1}
\end{align*}
and
\begin{align*}
 x= &X^{(1)}([b_{1}^{i}]_+)\Big(\sum\limits_{r=0}^{d_i}h_{i,r}(q^{\frac{1}{2}})q^{-\frac{1}{2}\Lambda_1(r[\beta_{1}^{i}]_+ +(d_i-r)[-\beta_{1}^{i}]_+,-e_i)} \\
 &\cdot q^{\frac{1}{2}\Lambda_1(-[b_{1}^{i}]_+,r[\beta_{1}^{i}]_+ +(d_i-r)[-\beta_{1}^{i}]_+)} X^{(1)}((r-d_i)\beta_{1}^{i})\Big)z^{-1}\\
=&q^{-\frac{1}{2}\Lambda_1([b_{1}^{i}]_{+},-e_{i})} X^{(1)}([b_{1}^{i}]_+)\Big(\sum\limits_{r=0}^{d_i}h_{i,r}(q^{\frac{1}{2}})q^{\frac{1}{2}(r-d_i)\widetilde{d}_id_{i}^{-1}} X^{(1)}((r-d_i)\beta_{1}^{i}) \Big)z^{-1}.
\end{align*}
It follows that
\begin{align}\label{L1}
L=& q^{\frac{1}{2}\Lambda_1([b_{1}^{i}]_{+},-e_{i})-\frac{1}{2}\Lambda_2([b_{2}^{i}]_{+},-e_i)+\frac{1}{2} \Lambda_1(E_{+}^{\widetilde{B}_1}[b_{2}^{i}]_{+},-[b_{1}^{i}]_{+})} X^{(1)}(E_{+}^{\widetilde{B}_1}[b_{2}^{i}]_{+}-[b_{1}^{i}]_{+}) \\
=&q^{\frac{1}{2}\Lambda_1([b_{1}^{i}]_{+},-e_{i})-\frac{1}{2}\Lambda_2([b_{2}^{i}]_{+},-e_i)+\frac{1}{2} \Lambda_1(E_{+}^{\widetilde{B}_1}[b_{2}^{i}]_{+},-[b_{1}^{i}]_{+})} X(E_{+}^{\widetilde{B}_1}[b_{2}^{i}]_{+}-[b_{1}^{i}]_{+}). \nonumber
\end{align}

When $b_{1,ij}<0$, by direct calculations, we have that
\begin{align*}
& y^{-1}Q(0)=X^{(1)}([b_{1}^{j}]_{+}-e_j), \\
& E_{-}^{\widetilde{B}_1}[-\beta_{2}^{i}]_+ =[-\beta_{2}^{i}]_+ +2\beta_{2,ji}e_j-\beta_{2,ji}[b_{1}^{j}]_+, \\
&  E_{-}^{\widetilde{B}_1}[\beta_{2}^{i}]_+ =[\beta_{2}^{i}]_+, \\
& E_{-}^{\widetilde{B}_1}\beta_{2}^{i}= \beta_{2}^{i} -2\beta_{2,ji}e_j+\beta_{2,ji}[b_{1}^{j}]_+ =\beta_{1}^{i}, \\
& E_{-}^{\widetilde{B}_1}e_j=[b_{1}^{j}]_{+} -e_j.
\end{align*}

It follows that
\begin{align*}
x=&X^{(1)}([-b_{1}^{i}]_{+})\Big(\sum\limits_{r=0}^{d_i}h_{i,r}(q^{\frac{1}{2}})q^{-\frac{1}{2}\Lambda_1 (r[\beta_{1}^{i}]_{+}+(d_i-r)[-\beta_{1}^{i}],-e_i)}\\
&\cdot q^{\frac{1}{2}\Lambda_1 (-[-b_{1}^{i}]_{+}, r[\beta_{1}^{i}]_{+} +(d_i-r)[-\beta_{1}^{i}]_{+})} X^{(1)}(r\beta_{1}^{i})\Big)z^{-1}\\
=&q^{-\frac{1}{2}\Lambda_1([-b_{1}^{i}]_{+},-e_i)}X^{(1)}([-b_{1}^{i}]_{+})\Big(\sum\limits_{r=0}^{d_i}h_{i,r} (q^{\frac{1}{2}})q^{\frac{1}{2}r\widetilde{d}_id_{i}^{-1}}X^{(1)}(r\beta_{1}^{i})\Big)z^{-1}
\end{align*}
and
\begin{align*}
&z^{-1}R(y^{-1}Q(0))\\
=&\sum\limits_{r=0}^{d_i}h_{i,r}(q^{\frac{1}{2}})q^{-\frac{1}{2}\Lambda_2(r[\beta_{2}^{i}]_{+} +(d_i-r)[-\beta_{2}^{i}]_{+}, -e_i)}\\
&\cdot q^{-\frac{1}{2}\Lambda_2(r[\beta_{2}^{i}]_{+} +(d_i-r)[-\beta_{2}^{i}]_{+} +(d_i-r)\beta_{2,ji}e_j, -(d_i-r)\beta_{2,ji}e_j)}\\
& \cdot q^{\frac{1}{2}\Lambda_1(r[\beta_{2}^{i}]_{+} +E_{-}^{\widetilde{B}_1}(d_i-r)[-\beta_{2}^{i}]_{+}, -E_{-}^{\widetilde{B}_1}(d_i-r)\beta_{2,ji}e_j)} \cdot X^{(1)}(r[\beta_{2}^{i}]_{+} +E_{-}^{\widetilde{B}_1}(d_i-r)[-\beta_{2}^{i}]_{+})z^{-1}\\
=&X^{(1)}(E_{-}^{\widetilde{B}_1}[-b_{2}^{i}]_{+}) \Big(\sum\limits_{r=0}^{d_i}h_{i,r}(q^{\frac{1}{2}}) q^{-\frac{1}{2}\Lambda_2(r[\beta_{2}^{i}]_{+} +(d_i-r)[-\beta_{2}^{i}]_{+}, -e_i)} X^{(1)}(r\beta_{1}^{i})\Big)z^{-1}\\
=&q^{-\frac{1}{2}\Lambda_2([-b_{2}^{i}]_{+}, -e_i)} X^{(1)}(E_{-}^{\widetilde{B}_1}[-b_{2}^{i}]_{+}) \Big( \sum\limits_{r=0}^{d_i}h_{i,r}(q^{\frac{1}{2}}) q^{\frac{1}{2}r\widetilde{d}_id_{i}^{-1}}X^{(1)}(r\beta_{1}^{i})\Big)z^{-1}.
\end{align*}
Then we have
\begin{align}\label{L2}
L =& q^{\frac{1}{2}\Lambda_2([-b_{2}^{i}]_{+}, e_i) -\frac{1}{2}\Lambda_1([-b_{1}^{i}]_{+},e_i) -\frac{1}{2}\Lambda_1(E_{-}^{\widetilde{B}_1}[-b_{2}^{i}]_{+}, [-b_{1}^{i}]_{+})}
 X^{(1)}(E_{-}^{\widetilde{B}_1}[-b_{2}^{i}]_{+}-[-b_{1}^{i}]_{+})  \\
= & q^{\frac{1}{2}\Lambda_2([-b_{2}^{i}]_{+}, e_i) -\frac{1}{2}\Lambda_1([-b_{1}^{i}]_{+},e_i) -\frac{1}{2}\Lambda_1(E_{-}^{\widetilde{B}_1}[-b_{2}^{i}]_{+}, [-b_{1}^{i}]_{+})}
 X(E_{-}^{\widetilde{B}_1}[-b_{2}^{i}]_{+}-[-b_{1}^{i}]_{+}). \nonumber
\end{align}
The proof is completed.
\end{proof}

\begin{lemma}\label{lemleadterm}
The cluster variables $u$ and $v$ belong to $\mathbb{L}_0$ with the leading terms being $y^{-1}Q(0)$ and $Lx$ respectively.
\end{lemma}
\begin{proof}
We analyze the following cases:

\begin{itemize}
  \item[(a)] If $b_{1,ij}=0$, then $u=y^{-1}Q(0)=\sum\limits_{r=0}^{d_j}h_{j,r}(q^{\frac{1}{2}})X(r[\beta^{j}_1]_+ +(d_j-r)[-\beta^j_1]_+ -e_j)$.   \item[(b)] If $b_{1,ij}>0$, then
\begin{align*}
u=&\sum\limits_{r=0}^{d_j}q^{-\frac{1}{2}\Lambda_1(r[\beta^{j}_1]_+ +(d_j-r)[-\beta^j_1]_+ -e_j,r\beta_{1,ij}e_i)} h_{j,r}(q^{\frac{1}{2}}) \\
&\cdot X(r[\beta^{j}_1]_+ +(d_j-r)[-\beta^j_1]_+ -e_j -r\beta_{1,ij}e_i)z^{r\beta_{1,ij}}
\end{align*}
and $y^{-1}Q(0)=X^{(1)}([-b^j_1]_+ -e_j)=X([-b^j_1]_+ -e_j)$.  \item[(c)] If $b_{1,ij}<0$, then
\begin{align*}
u=&\sum\limits_{r=0}^{d_j}q^{-\frac{1}{2}\Lambda_1(r[\beta^{j}_1]_+ +(d_j-r)[-\beta^j_1]_+ -e_j,(r-d_j)\beta_{1,ij}e_i)} h_{j,r}(q^{\frac{1}{2}}) \\
&\cdot X(r[\beta^{j}_1]_+ +(d_j-r)[-\beta^j_1]_+ -e_j +(d_j-r)\beta_{1,ij}e_i)z^{(r-d_j)\beta_{1,ij}}
\end{align*}
and $y^{-1}Q(0)=X([b^j_1]_+ -e_j)$. \end{itemize}
Note that $z=X^{(1)}_i\in \mathbb{L}_0$. It follows that $u\in\mathbb{L}_0$.

Therefore $y^{-1}Q(x^{-1}P(y))-y^{-1}Q(0)\in\bigoplus\limits_{l>0}^{\infty}\mathbb{L}_{0}^{\vee}x^{-l}$ and $y^{-1}Q(0)\in\mathbb{L}_{0}^{\vee}$. Hence $y^{-1}Q(0)$ is the leading term of $u$.

Note that
\begin{align*}
v=&z^{-1}R(y^{-1}Q(z))=z^{-1}(R(y^{-1}Q(z))-R(y^{-1}Q(0)))+z^{-1}R(y^{-1}Q(0))\\
=&z^{-1}(R(y^{-1}Q(z))-R(y^{-1}Q(0)))+Lx.
\end{align*}

By Lemma \ref{lemL}, we know that $Lx\in\mathbb{L}_{0}^{\vee}x$. Note that $z=x^{-1}P(y)\in\mathbb{L}_{0}^{\vee}x^{-1}$ and
\begin{align*}
R(y^{-1}Q(z))-R(y^{-1}Q(0))\in \bigoplus\limits_{l>0}^{\infty}\mathbb{L}_{0}^{\vee}z^{l},
\end{align*}
which implies that
\begin{align*}
z^{-1}(R(y^{-1}Q(z))-R(y^{-1}Q(0)))\in\bigoplus\limits_{l\geq0}^{\infty}\mathbb{L}_{0}^{\vee}z^{l} \subseteq \bigoplus\limits_{l\leq0}^{\infty}\mathbb{L}_{0}^{\vee}x^{l} =\bigoplus\limits_{l<1}^{\infty}\mathbb{L}_{0}^{\vee}x^{l}.
\end{align*}
Hence $v\in\mathbb{L}_0$ and $Lx$ is the leading term of $v$.
\end{proof}

In the following, we need the divisibility in the noncommutative domain $A$. For $a,b\in A$, we say that $a$ is left divided (resp. right divided) by $b$ if there exits an element $c\in A$ (resp. $d\in A$) such that $bc=a$ (resp. $db=a$), which is denoted by $b|_{l}a$ (resp. $b|_ra$).

\begin{lemma}\label{mainlem}
If $z^{N}F=Gu^{N_1}v^{N_2}$ for $F,G\in\mathbb{L}_0\setminus\{0\}$ and positive integers $N,N_1$ and $N_2$, then $G$ is left divided by $z^N$ in $\ll_0$.
\end{lemma}
\begin{proof}
Let $\sigma$ be an automorphism of $\A$ defined by $\sigma(X_k)=q^{\lambda_{0,ik}}X_k$ for all $k\in[1,m]$. It follows that $\sigma(X(c))=q^{\Lambda_0(e_i,c)}X(c)$ and $xX(c)=\sigma(X(c))x$ for $c\in\ZZ^{m}$. We obtain that
$$
z^N=(x^{-1}P)^N=\sigma^{-1}(P)\sigma^{-2}(P)\cdots\sigma^{-N}(P)x^{-N}.
$$
By Lemma \ref{lemleadterm}, it follows that
$$
\sigma^{-1}(P)\sigma^{-2}(P)\cdots\sigma^{-N}(P)x^{-N}F=G\big((y^{-1}Q(0))^{N_1}(Lx)^{N_2} +\rm{lower~order~terms}\big).
$$
There exists some $F^\prime\in\mathbb{L}_0\setminus\{0\}$ such that
$$
\sigma^{-1}(P)\cdots\sigma^{-N}(P)F^{\prime}=G\big((y^{-1}Q(0))^{N_1}x^{N_2} +\rm{lower~order~terms}\big).
$$

We assume that $G$ cannot be left divided by $z^N$. The element $G$ can be written as $G=\sum\limits_{l\in\ZZ}g_lx^l$, where $g_l\in\ll_{0}^{\vee}$. We can assume that $l^{\prime}$ is the maximal integer such that $g_{l^\prime}$ is not left divided by $\sigma^{-1}(P)\cdots\sigma^{-N}(P)$ in $\ll_{0}^{\vee}$. Note that $g_{l^\prime}\sigma^{l^{\prime}}((y^{-1}Q(0))^{N_1})x^{l^{\prime}+N_2}$ is left divided by $\sigma^{-1}(P)\cdots\sigma^{-N}(P)$ in $\ll_{0}$. Hence $g_{l^\prime}\sigma^{l^{\prime}}((y^{-1}Q(0))^{N_1})$ is left divided by $\sigma^{-1}(P)\cdots\sigma^{-N}(P)$ in $\ll_{0}^{\vee}$ since $\sigma^{-1}(P)\cdots\sigma^{-N}(P)\in\ll_{0}^{\vee}$.

Recall that $y^{-1}Q(0)=X([-b_{1}^{j}]_{+}-e_j)$ if $b_{1,ij}>0$ and $y^{-1}Q(0)=X([b_{1}^{j}]_{+}-e_j)$ if $b_{1,ij}<0$.
If $b_{1,ij}\neq0$ then $y^{-1}Q(0)$ is invertible in $\ll_0^{\vee}$, which implies that $g_{l^{\prime}}$ is left divided by  $\sigma^{-1}(P)\cdots\sigma^{-N}(P)$ in $\ll_{0}^{\vee}$. It is a contradiction.

In the rest of the proof, assume that $b_{1,ij}=0$. Note that $Q(0)=Q$ and $\beta_{1}^{j}=\beta_{0}^{j}$.
Using the fact that the principal part of $\widetilde{B}_1$ is skew-symmetrizable, it follows that $b_{1}^{i}, b_{1}^{j}\in\sum\limits_{k\in[1,m] \setminus\{i,j\}}\ZZ e_k$. Hence $\beta_{1}^{i}, \beta_{1}^{j}\in\sum\limits_{k\in[1,m] \setminus\{i,j\}}\ZZ e_k$.
By (\ref{equ1}), it follows that $\sigma(X(\beta_{0}^{j}))=X(\beta_{0}^{j})$ and
$\sigma(X(\beta_{0}^{i}))=q^{-\widetilde{d}_id^{-1}_i}X(\beta_{0}^{i})$.

Recall that $\beta_{1}^{i} =-\beta_{0}^{i}$. Since $b_{1,ij}=0$, it follows that $\beta^{j}_1=\beta^j_0$,
$$
y^{-1}Q(0)=\sum\limits_{r=0}^{d_j}h_{j,r}(q^{\frac{1}{2}})q^{\frac{1}{2}r\widetilde{d}_{j}d_{j}^{-1}}X(r\beta_{0}^{j}) X([-b_{0}^{j}]_+ -e_j),
$$
and
$$
P=\sum\limits_{r=0}^{d_i}h_{i,r}(q^{\frac{1}{2}})q^{\frac{1}{2}\Lambda_0 (e_i,[-b_{0}^{i}]_+) -\frac{1}{2}r\widetilde{d}_id_{i}^{-1}}X(r\beta_{0}^{i})X([-b_{0}^{i}]_+).
$$
By direct calculations, we conclude that
\begin{align*}
& \sigma^{l^\prime}((y^{-1}Q(0))^{N_1}) =\Big( \prod\limits_{k=1}^{N_1}\big(\sum\limits_{r=0}^{d_j} h_{j,r} (q^{\frac{1}{2}}) q^{\frac{2k-1}{2}r\widetilde{d}_jd^{-1}_j}X(r\beta_{0}^{j})\big)\Big) q^{-\frac{1}{2}N_{1}^{2}\widetilde{d}_j}\alpha \\
= & \Big( \prod\limits_{k=1}^{N_1}\big(\sum\limits_{r=0}^{d_j} h_{j,r} (q^{\frac{1}{2}}) q^{(k-\frac{1}{2})(rd_{j}^{-1}-1)\widetilde{d}_j}X(r\beta_{0}^{j})\big)\Big) \alpha
\end{align*}
and
\begin{align*}
& \sigma^{-1}(P)\cdots\sigma^{-N}(P)
=\Big( \prod_{k=1}^{N}\big(\sum\limits_{r=0}^{d_i} h_{i,r}(q^{\frac{1}{2}}) q^{\frac{2k-1}{2}r\widetilde{d}_{i}d^{-1}_i} X(r\beta_{0}^{i})\big)\Big) q^{-\frac{1}{2}N^2\widetilde{d}_i}\beta \\
= &\Big( \prod_{k=1}^{N}\big(\sum\limits_{r=0}^{d_i} h_{i,r}(q^{\frac{1}{2}}) q^{(k-\frac{1}{2})(rd_{i}^{-1}-1)\widetilde{d}_{i}} X(r\beta_{0}^{i})\big)\Big) \beta,
\end{align*}
where
$$
\alpha=q^{\frac{1}{2}N_{1}^{2}\widetilde{d}_j+l^{\prime}N_{1}\Lambda_{0} (e_i,[-b_{0}^{j}]_+-e_j)}  X([-b_{0}^{j}]_+ -e_j)^{N_1}
$$
and
$$
\beta=q^{\frac{1}{2}N^2\widetilde{d}_i -\frac{1}{2}N^2\Lambda_{0}(e_i,[-b_{0}^{i}]_+)} X([-b_{0}^{i}]_+)^{N}
$$
are Laurent monomials in $X_k$ for $k\in[1,m]\setminus\{i\}$. It is clear that $\alpha$ and $\beta$ are invertible in $\ll_{0}^{\vee}$. Let
\begin{align*}
V_1(t)=\prod\limits_{k=1}^{N_1}\big(\sum\limits_{r=0}^{d_j} h_{j,r} (q^{\frac{1}{2}}) q^{(k-\frac{1}{2})(rd_{j}^{-1}-1)\widetilde{d}_j}t^r\big),
\end{align*}
\begin{align*}
V_2(t)= \prod_{k=1}^{N}\big(\sum\limits_{r=0}^{d_i} h_{i,r}(q^{\frac{1}{2}}) q^{(k-\frac{1}{2})(rd_{i}^{-1}-1)\widetilde{d}_{i}} t^r\big) \in\ZZ[q^{\pm\frac{1}{2}}][t].
\end{align*}
Then $V_2(X(\beta_{0}^{i})) |_{l} g_{l^{\prime}} V_1(X(\beta_{0}^{j}))$ in $\ll_{0}^{\vee}$.

By (\ref{equ1}), we know that $\Lambda_0(\beta_{0}^{i},e_k)=\Lambda_0(\beta_{0}^{j},e_k)=0$ for $k\in[1,m]\setminus\{i,j\}$. Hence $X(\beta_{0}^i)$ and $X(\beta_{0}^{j})$ are in the center of $\ll^{\vee\vee}_0$.

Let $\varphi$ be an automorphism of $\ll_{0}^{\vee\vee}$ defined by $\varphi(X_k)=q^{\lambda_{0,jk}}X_k$ for $k\in[1,m]\setminus\{i,j\}$.
Then we obtain that
$$
yX(\beta^{j}_0)=\varphi(X(\beta^{j}_0))y=q^{\Lambda_0(e_j,\beta^{j}_0)}X(\beta^{j}_0)y =q^{-\widetilde{d}_jd_{j}^{-1}}X(\beta^{j}_0)y.
$$
Let $g_{l^{\prime}}=\sum\limits_{r\in\ZZ}g_{l^{\prime},r}y^{r}$, where $g_{l^{\prime},r}\in\ll_{0}^{\vee\vee}$ for all $r\in\ZZ$. Note that $V_2(X(\beta_{0}^{i}))\in \ll^{\vee\vee}_{0}$.
Thus $V_2(X(\beta_{0}^{i})) |_{l} g_{l^{\prime},r} V_1(q^{-r\widetilde{d}_jd_{j}^{-1}} X(\beta_{0}^{j}))$ in $\ll_{0}^{\vee\vee}$ for all $r\in\ZZ$.

Because the matrix $\widetilde{B}_0$ is of full rank, we have that $\beta_{0}^{i}$ and $\beta_{0}^{j}$ are $\ZZ$-linearly independent. Using Lemma \ref{coprimelem}, it implies that $V_1(X(\beta_{0}^{i}))$ and $V_2(q^{-r\widetilde{d}_jd_{j}^{-1}} X(\beta_{0}^{j}))$ are coprime in the center of $\ll_{0}^{\vee\vee}$.
It follows that $V_2(X(\beta_{0}^{i})) |_{l}  \,  g_{l^{\prime},r}$ in $\ll_{0}^{\vee\vee}$ and $V_2(X(\beta_{0}^{i})) |_{l} \, g_{l^{\prime}}$ in $\ll_{0}^{\vee}$. Thus $g_{l^{\prime}}$ is left divided by $\sigma^{-1}(P)\vdots\sigma^{-1}(P)$ in $\ll_{0}^{\vee}$, which is a contradiction. The proof is completed.
\end{proof}

\begin{proof}[{\bf{Proof of Theorem \ref{laurentphenomenon}}}]
It suffices to prove that any cluster variable $Z$ belongs to $\ll_0$. By induction hypothesis, we can assume that $Z$ belongs to $\ll_1$ and $\ll_3$. Since $Z\in\ll_1$ and $\ll_3$, we have that $Z=z^{-N}G=Fu^{-N_1}v^{-N_2}$ for $F,~G\in\ll_0$ and $N,~N_1,~N_2\geq0$. By Lemma \ref{mainlem}, we have that $G$ is left divided by $z^{N}$. It follows that $Z=z^{-N}G\in\ll_0$.
\end{proof}

\begin{example}(Type $\emph{G}_2$)\label{example1}
Let $\A(1,3)$ denote the generalized quantum cluster algebra associated with the initial seed $(X,\mathbf{h},\Lambda,B)$, where $\mathbf{d}=(3,1)$,
\begin{equation*}
\Lambda=\left(
  \begin{array}{cc}
    0 & 1 \\
    -1 & 0 \\
  \end{array}
\right)
\text{~and~}
B=\left(
  \begin{array}{cc}
    0 & 1 \\
    -3 & 0 \\
  \end{array}
\right).
\end{equation*}
By the definition, $\A(1,3)$ is the $\ZZ[q^{\pm\frac{1}{2}}]$-subalgebra of $\mathcal{F}$ generated by the cluster variables $\{X_i~|~i\in\ZZ\}$ which are obtained from the exchange relations:
\begin{equation*}
   X_{m-1}X_{m+1}=\begin{cases}
                   q^{\frac{1}{2}}X_m+1, & \text{ if $m$ is odd,}\\
                   q^{\frac{3}{2}}X^{3}_{m} +h(q^{\frac{1}{2}})qX^{2}_{m} +h(q^{\frac{1}{2}})q^{\frac{1}{2}}X^{}_{m} +1, & \text{ if $m$ is even,}
                  \end{cases}
\end{equation*}
for $h(q^{\frac{1}{2}})\in\ZZ[q^{\pm\frac{1}{2}}]$. It is clear that $X_mX_{m+1}=qX_{m+1}X_m$ for $m\in\ZZ$.

Recall that $X(a,b)=q^{-\frac{ab}{2}}X_{1}^{a}X_{2}^{b}$ for any $a,b\in\ZZ$. By direct calculations, it follows that
\begin{enumerate}
\item[] $
X_3=X(-1,0)+h(q^{\frac{1}{2}})X(-1,1)+h(q^{\frac{1}{2}})X(-1,2)+X(-1,3);
$
\item[] $
X_4=X(-1,-1)+h(q^{\frac{1}{2}})X(-1,0)+h(q^{\frac{1}{2}})X(-1,1)+X(-1,2)+X(0,-1);
$
\item[] $
X_5=X(-2,-3)+h(q^{\frac{1}{2}})(q^{-\frac{1}{2}}+q^{\frac{1}{2}})X(-2,-2) \\ {\hspace{1.0cm}}+[h(q^{\frac{1}{2}})(q^{-1}+q)+h(q^{\frac{1}{2}})^2]X(-2,-1) \\ {\hspace{1.0cm}}+[q^{-\frac{3}{2}}+q^{\frac{3}{2}}+h(q^{\frac{1}{2}})^2(q^{-\frac{1}{2}}+q^{\frac{1}{2}})]X(-2,0) \\ {\hspace{1.0cm}}+[h(q^{\frac{1}{2}})(q^{-1}+q^{})+h(q^{\frac{1}{2}})^2]X(-2,1) +h(q^{\frac{1}{2}})(q^{-\frac{1}{2}}+q^{\frac{1}{2}})X(-2,2) \\
{\hspace{1.0cm}}+X(-2,3) +(q^{-1}+1+q)X(-1,-3) +h(q^{\frac{1}{2}})(q^{-1}+2+q)X(-1,-2)\\
{\hspace{1.0cm}}+h(q^{\frac{1}{2}})(q^{-1}+1+q+h(q^{\frac{1}{2}}))X(-1,-1) +(q^{-1}+1+q+h(q^{\frac{1}{2}})^2)X(-1,0)\\
{\hspace{1.0cm}}+h(q^{\frac{1}{2}})X(-1,1) +(q^{-1}+1+q)X(0,-3)\\
{\hspace{1.0cm}}+h(q^{\frac{1}{2}})(q^{-\frac{1}{2}} +q^{\frac{1}{2}})X(0,-2) +h(q^{\frac{1}{2}})X(0,-1)+X(1,-3);
$
\item[] $
X_6=X(-1,-2) +h(q^{\frac{1}{2}})X(-1,-1) +h(q^{\frac{1}{2}})X(-1,0) +X(-1,1)\\
{\hspace{1.0cm}} +(q^{-\frac{1}{2}}+q^{\frac{1}{2}})X(0,-2) +h(q^{\frac{1}{2}})X(0,-1) +X(1,-2);
$
\item[] $
X_7=X(-1,-3)+h(q^{\frac{1}{2}})X(-1,-2) +h(q^{\frac{1}{2}})X(-1,-1) +X(-1,0)\\
{\hspace{1.0cm}} +(q^{-1}+1+q)X(0,-3) +h(q^{\frac{1}{2}})(q^{-\frac{1}{2}}+q^{\frac{1}{2}})X(0,-2) +h(q^{\frac{1}{2}})X(0,-1) \\
{\hspace{1.0cm}} +(q^{-1}+1+q)X(1,-3) +h(q^{\frac{1}{2}})X(1,-2) +X(2,-3);
$
\item[] $
X_8=X_0=X(0,-1)+X(1,-1);
$
\item[] $
X_9=X_1;
$
\item[] $
X_{10}=X_2.
$
\end{enumerate}

\end{example}

Example \ref{example1} gives the explicit cluster expansion formulas of every cluster variable associated to the initial  cluster for type $\emph{G}_2$. Note that all coefficients in these cluster expansion formulas belong to $\NN[q^{\pm\frac{1}{2}},h(q^{\frac{1}{2}})]$ which lead to the following conjecture:
\begin{conjecture}[Positivity conjecture]
Any cluster variable $X$ of  a generalized quantum cluster algebra $\mathcal{A}(X, \textbf{h}, \Lambda, \widetilde{B})$ can be expressed as a Laurent polynomial $f(X_{1}^{\pm1},\ldots,X_{m}^{\pm1})$ in the initial extended cluster $\{X_{1}^{\pm1},\ldots,X_{m}^{\pm1}\}$ with the coefficients belonging to $\NN[q^{\pm\frac{1}{2}},h_{i,r}(q^{\frac{1}{2}})]$.
\end{conjecture}

\section{Upper bounds}

Let $(X, \textbf{h}, \Lambda, \widetilde{B})$ be a quantum seed in the skew-field $\mathcal{F}$. For simplicity, we set  $\mathbb{L}(X, \textbf{h}, \Lambda, \widetilde{B})=\ZZ\mathbb{P}[X^{\pm1}_1,\ldots,,X^{\pm1}_{n}]$. As in \cite[Definition 1.1]{bfz}, the upper bound is defined by
$$
\mathcal{U}(X, \textbf{h}, \Lambda, \widetilde{B}):=\mathbb{L}(X, \textbf{h}, \Lambda, \widetilde{B})\cap \bigcap\limits_{i=1}^{n}\mathbb{L}(\mu_i(X, \textbf{h}, \Lambda, \widetilde{B})).
$$

For $i\neq j$, two polynomials $X_iX^\prime_i$ and $X_jX^\prime_j$ are called coprime if there does not exist any non-invertible element $c$ in the center of $\ZZ\PP[X_1,\ldots,X_n]$ such that both $X_iX^\prime_i$ and $X_jX^\prime_j$ are divided by $c$. The quantum seed $(X,\textbf{h},\Lambda,\widetilde{B})$ is said to be coprime if $X_iX^\prime_i$ and $X_jX^\prime_j$ are coprime for $i,j\in[1,n]$ and $i\neq j$.

Let $V_{b^{i}}^0=1$ and $V_{b^{i}}^s = \prod\limits_{k=1}^{s} (\sum\limits_{l=0}^{d_i} q^{\frac{2k-1}{2}l\widetilde{d}_id^{-1}_i} h_{i,l}(q^{\frac{1}{2}})X(l\beta^{i}) )$, where $i\in[1,n]$ and $s\in\ZZ_+$.
Similarly, we define $W_{b^{i}}^0=V_{b^{i}}^0=1$ and $W_{b^{i}}^s = \prod\limits_{k=1}^{s} (\sum\limits_{l=0}^{d_i} q^{\frac{1-2k}{2}l\widetilde{d}_id^{-1}_i} h_{i,l}(q^{\frac{1}{2}})X(-l\beta^{i}) )$ for $i\in[1,n]$ and $s\in\ZZ_+$.
\begin{lemma}\label{lemxprime}
For $i\in[1,n]$ and $s\in\ZZ_+$, we have that
\begin{equation}\label{equxprime}
(X^{\prime}_i)^s=V_{b^{i}}^s X(s[-b^i]_+-se_i),
\end{equation}
and
\begin{equation}\label{equxprime2}
(X^{\prime}_i)^s=W_{b^{i}}^s X(s[b^i]_+-se_i).
\end{equation}
\end{lemma}
\begin{proof}
Since $\Lambda([-b^{i}]_+ -e_i,r\beta^{i})=\Lambda_0(-e_i,r\beta^{i})=r\widetilde{d}_id^{-1}_i$, it follows that
$$
X^{\prime}_i=\left(\sum\limits_{r=0}^{d_i}h_{i,r}(q^{\frac{1}{2}})q^{\frac{1}{2}r\widetilde{d}_id^{-1}_i}X(r\beta^i)\right) X([-b^i]_+-e_i)
$$
and
$$
X([-b^i]_i-e_i)X(r\beta^i)=q^{r\widetilde{d}_id^{-1}_i}X(r\beta^i)X([-b^i]_+-e_i).
$$
Then Equality (\ref{equxprime}) holds.

Since $h_{i,l}(q^{\frac{1}{2}})=h_{i,d_i-l}(q^{\frac{1}{2}})$ for each $l\in[0,d_i]$, we can write
$$
X^\prime_i=\sum\limits_{r=0}^{d_i}h_{i,l}(q^{\frac{1}{2}})X(r[-\beta^i]_++(d_i-r)[\beta^i]_+-e_i).
$$
The rest of the proof of Equality (\ref{equxprime2}) runs similarly as Equality~(\ref{equxprime}).
\end{proof}

The following result is a generalization of \cite[Proposition 5.4]{BZ2005}.
\begin{proposition}\label{propu1}
We have
$$
\mathcal{U}(X, \textbf{h}, \Lambda, \widetilde{B})=\bigcap\limits_{i=1}^{n}\ZZ\PP[X^{\pm1}_1,\ldots,X^{\pm1}_{i-1},X_i,X^{\prime}_i, X^{\pm1}_{i+1},\ldots,X^{\pm1}_{n}].
$$
\end{proposition}
\begin{proof}
We only need to show that
$$
\ZZ\PP[X^{\pm1}_1,\ldots,X^{\pm1}_n]\cap\ZZ\PP[(X^\prime_1)^{\pm1},X^{\pm1}_2,\ldots,X^{\pm1}_n] =\ZZ\PP[X_1,X^\prime_1,,X^{\pm1}_2,\ldots,X^{\pm1}_n].
$$
It suffices to prove the ``$\subseteq$" inclusion.

Let $Y\in\ZZ\PP[X^{\pm1}_1,X^{\pm1}_2,\ldots,X^{\pm1}_n]\cap\ZZ\PP[(X^{\prime}_1)^{\pm1},X^{\pm1}_2,\ldots,X^{\pm1}_n]$. We claim that
\begin{equation}\label{equy}
Y=c_0+\sum\limits_{i\in\ZZ_+}(c_iX^{i}_1+c^{\prime}_i(X^\prime_1)^{i}),
\end{equation}
where $c_0,c_i,c^{\prime}_i\in\ZZ\PP[X^{\pm1}_2,\ldots,X^{\pm1}_n]$ and all but finitely many $c_i,c^{\prime}_i$ are zero.
Using the relation (\ref{qcommutation}) if necessary, we obtain that $Y=\sum\limits_{i\in\ZZ}c_iX^{i}_1$, where $c_i\in\ZZ\PP[X^{\pm1}_2,\ldots,X^{\pm1}_n]$ and $c_i=0$ for all but finitely many $i\in\ZZ$.
Assume that there exists some integer $j>0$ such that $c_{-j}\neq0$. By Lemma \ref{lemxprime}, we have that
$$
(X_1)^{-j}=q^{-\frac{1}{2}\Lambda(j[-b^1]_+,je_1)}(V^{j}_{b^1}X(j[-b^1]_+))^{-1}(X^\prime_1)^j.
$$
Since $Y\in\ZZ\PP[(X^\prime_1)^{\pm1},X^{\pm1}_2,\ldots,X^{\pm1}_n]$, $c_{-j}$ is right divided by $V^{j}_{b^1}X(j[-b^1]_+)$. Let $c^\prime_j=q^{-\frac{1}{2}\Lambda(j[-b^1]_+,je_1)}c_{-j}(V^{j}_{b^1}X(j[-b^1]_+))^{-1}$. Hence (\ref{equy}) holds. Then we conclude that $Y\in\ZZ\PP[X_1,X^\prime_1,,X^{\pm1}_2,\ldots,X^{\pm1}_n]$. The proof is completed.
\end{proof}

From now on, let $(X,\mathbf{h},\Lambda_0,\widetilde{B}_0)$ be the initial quantum seed. For simplicity, define $(X^{\prime},\mathbf{h},\Lambda_1,\widetilde{B}_1)=\mu_i(X,\mathbf{h},\Lambda_0,\widetilde{B}_0)$ and $(X^{\prime\prime},\mathbf{h},\Lambda_2,\widetilde{B}_2)=\mu_j(X^{\prime},\mathbf{h},\Lambda_1,\widetilde{B}_1)$, where $i,j\in[1,n]$ and $i\neq j$. Recall that $X^{\prime\prime}_j=\mu_j\mu_i(X_j)$.

In the next, from  Lemma \ref{lemcap} to Lemma \ref{lemcapn=2}, we only consider the case when $n=2$.

\begin{lemma}\label{lemcap}
We have
$
\ZZ\PP[X_1,X^{\pm1}_2]\cap\ZZ\PP[X^{\pm1}_1,X_2,X^\prime_2]=\ZZ\PP[X_1,X_2,X^\prime_2].
$
\end{lemma}
\begin{proof}
Obviously, $\ZZ\PP[X_1,X^{\pm1}_2]\cap\ZZ\PP[X^{\pm1}_1,X_2,X^\prime_2]\supseteq\ZZ\PP[X_1,X_2,X^\prime_2]$.

To show the other inclusion, let $Y\in \ZZ\PP[X_1,X^{\pm1}_2]\cap\ZZ\PP[X^{\pm1}_1,X_2,X^\prime_2]$. It is easily seen that $Y=\sum\limits_{m\in\ZZ}c_mX^{m}_1$ for $c_m\in\ZZ\PP[X_2,X^\prime_2]$. Note that $c_m=0$ for all but finitely many $m\in\ZZ$.

If $c_m\neq0$, then we can write $c_m=c_{m,0}+c_{m,1}(X_2)+c_{m.2}(X^\prime_2)$, where $c_{m,0}\in\ZZ\PP$, $c_{m,1}$ and $c_{m.2}$ are polynomial over $\ZZ\PP$ which have the constant term $0$. Let $m_0$ be the minimal integer which satisfies the condition $c_{m_0}\neq0$. If $m_0\geq0$, then the statement is true. If $m_0<0$, then
$$
c_{m_0}=c_{m_0,0}+c_{m_0,1}(X_2)+c_{m_0,2}\big(\sum\limits_{r=0}^{d_2}h_{2,r}(q^{\frac{1}{2}}) X(r[\beta^2_0]_++(d_2-r)[-\beta^2_0]_+-e_2)\big).
$$
Hence $0\neq c_{m_0}X^{m_0}_1\notin\ZZ\PP[X_1,X^{\pm1}_2]$ and $Y\notin\ZZ\PP[X_1,X^{\pm1}_2]$, a contradiction.
\end{proof}

\begin{lemma}\label{lemplus}
If $b_{0,12}\neq0$, then $\ZZ\PP[X_1,X^\prime_1,X^{\pm1}_2]=\ZZ\PP[X_1,X^\prime_1,X_2,X^\prime_2]+\ZZ\PP[X_1,X^{\pm1}_2]$.
\end{lemma}
\begin{proof}
It is clear that $\ZZ\PP[X_1,X^\prime_1,X^{\pm1}_2]\supseteq\ZZ\PP[X_1,X^\prime_1,X_2,X^\prime_2]+\ZZ\PP[X_1,X^{\pm1}_2]$.

Using (\ref{qcommutation}) if necessary, it is clear that there exist $Y_1,Y_2\in\ZZ\PP[X_1,X^\prime_1,X^{\pm1}_2]$ such that $YX^{-1}_2=X^{-1}_2Y_1$ and $X^{-1}_2Y=Y_2X^{-1}_2$ for any $Y\in\ZZ\PP[X_1,X^\prime_1,X^{\pm1}_2]$.
To prove the reverse inclustion, it suffices to show that
$$
(X^\prime_1)^kX^{-l}_2,X^{-l}_2(X^\prime_1)^k\in\ZZ\PP[X_1,X^\prime_1,X_2,X^\prime_2]+\ZZ\PP[X_1,X^{\pm1}_2]
$$
for all positive integers $k$ and $l$.

Case 1. If $b_{0,12}>0$, then $b_{0,21}<0$. Since $X^\prime_2=\sum\limits_{r=0}^{d_2}h_{2,r}(q^{\frac{1}{2}})X(r\beta_0^{2} +([-b_0^2]_+-e_2)$,
we have that
$$
X^{-1}_2+\sum\limits_{r=1}^{d_2}h_{2,r}(q^{\frac{1}{2}})q^{\frac{1}{2}r\widetilde{d}_2d^{-1}_{2}}X(r\beta^{2}_0)X^{-1}_2 =q^{-\frac{1}{2}\Lambda_0([-b^2_0]_+,e_2)}X(-[-b^2_0]_+)X^{\prime}_2.
$$
We set $P=-\sum\limits_{r=1}^{d_2}h_{2,r}(q^{\frac{1}{2}})q^{\frac{1}{2}r\widetilde{d}_2d^{-1}_{2}}X(r\beta^{2}_0)$ and $R=q^{-\frac{1}{2}\Lambda_0([-b^2_0]_+,e_2)}X(-[-b^2_0]_+)$. Note that $P\in\ZZ\PP[X_1]$ and $R,R^{-1}\in\ZZ\PP$. Then
\begin{align*}
&X^{-1}_2=PX^{-1}_2+RX^\prime_2=P^2X^{-1}_2+PRX^{\prime}_2+RX^\prime_2 \\
=&P^3X^{-1}_2+P^2RX^\prime_2+PRX^\prime_2+RX^\prime_2=\ldots \\
=&P^kX^{-1}_2+(P^{k-1}+P^{k-2}+\ldots+P+1)RX^\prime_2
=P^kX^{-1}_2+Q.
\end{align*}

Set $Q=(P^{k-1}+P^{k-2}+\ldots+P+1)RX^\prime_2$. Then $Q\in\ZZ\PP[X_1,X^\prime_2]$. It is easy to see that $P$ is left and right divided by $X^{\beta_{0,12}}_1$, we obtain that $P^kX^{-1}_2$ is left and right divided by $X^{k\beta_{0,12}}_1$ by (\ref{qcommutation}). It follows that $P^kX^{-1}_2Q$ and $QP^kX^{-1}_2$ are left and right divied by $X^{k\beta_{0,12}}_1$, which implies that $(X^{-l}_2-Q^l)$ is left and right divided by $X^{k\beta_{0,12}}_1$. Hence we conclude that $(X^\prime_1)^{k}(X^{-l}_2-Q^l),(X^{-l}_2-Q^l)(X^\prime_1)^{k}\in\ZZ\PP[X_1,X^{\pm1}_2]$.

Because $Q^l\in\ZZ\PP[X_1,X^\prime_2]$, we have that $(X^\prime_1)^{k}Q^l,Q^l(X^\prime_1)^{k}\in\ZZ\PP[X_1,X^{\prime}_1,X_2,X^\prime_2]$. Therefore $(X^\prime_1)^{k}X^{-l}_2,X^{-l}_2(X^\prime_1)^{k}\in \ZZ\PP[X_1,X^{\pm1}_2]+\ZZ\PP[X_1,X^{\prime}_1,X_2,X^\prime_2]$, as desired.

Case 2. When $b_{0,12}<0$, write $X^\prime_2=\sum\limits_{r=0}^{d_2}h_{2,r}(q^{\frac{1}{2}})X(-r\beta_0^{2} +([b_0^2]_+-e_2)$. Then
$$
X^{-1}_2=-\sum\limits_{r=1}^{d_2}h_{2,r}(q^{\frac{1}{2}})q^{-\frac{1}{2}r\widetilde{d}_2d^{-1}_2}X(-r\beta^{2}_0)X^{-1}_2 +q^{-\frac{1}{2}\Lambda_0([b^{2}_0]_+ ,e_2)}X(-[b^2_0]_+)X^\prime_2.
$$
The remainder of the proof is similar to the one given above for Case 1.
\end{proof}

The following lemma is an extension of \cite[(4.4)]{bfz} and \cite[(5.15)]{BZ2005}.
\begin{lemma}\label{lemcapn=2}
If $(X,\mathbf{h},\Lambda_0,\widetilde{B}_0)$ is coprime, then
$$
\ZZ\PP[X_1,X^\prime_1,X_2,X^\prime_2]=\ZZ\PP[X_1,X^\prime_1,X^{\pm1}_2]\cap\ZZ\PP[X^{\pm1}_1,X_2,X^\prime_2].
$$
\end{lemma}
\begin{proof}
The lemma may be proved in much the same way as the proof of \cite[(4.4)]{bfz}. For readers' convenience, we sketch the proof here.

It suffices to show $\ZZ\PP[X_1,X^\prime_1,X_2,X^\prime_2]\supseteq\ZZ\PP[X_1,X^\prime_1,X^{\pm1}_2]\cap\ZZ\PP[X^{\pm1}_1,X_2,X^\prime_2]$. Let $Y\in\ZZ\PP[X_1,X^\prime_1,X^{\pm1}_2]\cap\ZZ\PP[X^{\pm1}_1,X_2,X^\prime_2]$.

Case 1. If $b_{0,12}=0$, then $b_{0,21}=0$, $V^{s}_{b^{1}_0}X_2=X_2V^{s}_{b^{1}_0}$ and $V^{t}_{b^{2}_0}X_1=X_1V^{t}_{b^{2}_0}$. Note that $V^{s}_{b^{1}_0}$ and $V^{t}_{b^{2}_0}$ are in the center of $\ZZ\PP$ for positive integers $s$ and $t$.

Let $Y=\sum\limits_{m_1,m_2\in\ZZ}c_{m_1,m_2}X^{m_1}_1X^{m_2}_2$, where all but finitely many $c_{m_1,m_2}\in\ZZ\PP$ are zero. Note that $(X^\prime_i)^s=q^{-\frac{1}{2}\Lambda_0(s[-b^i_0]_+,-se_i)}V^{s}_{b^i_0}X(s[-b^i_0]_+)X(-se_i)$ for $i=1,2$, and $q^{-\frac{1}{2}\Lambda_0(s[-b^i_0]_+,-se_i)}V^{s}_{b^i_0}X(s[-b^i_0]_+)\in\ZZ\PP$.

We only need to consider the case when $m_1<0$ and $m_2<0$.
When $m_1,m_2<0$, since $V^{m_1}_{b^1_0}$ and $V^{m_2}_{b^2_0}$ are coprime, we obtain that $c_{m_1,m_2}$ is divided by $V^{m_1}_{b^1_0}V^{m_1}_{b^2_0}$ in $\ZZ\PP$ by \cite[Proposition A.2]{BZ2005}. Therefore $Y\in\ZZ\PP[X_1,X^\prime_1,X_2,X^\prime_2]$, as desired.

Case 2. Assume that $b_{0,12}\neq0$. By Lemmas \ref{lemcap} and \ref{lemplus}, the rest of proof is similar to the proof of \cite[Lemma 4.5]{bfz} and so is omitted.
\end{proof}

The following result is a generalization of \cite[Proposition 5.7]{BZ2005}.
\begin{proposition}\label{propu2}
If $(X,\mathbf{h},\Lambda,\widetilde{B})$ is coprime, then
$$
\mathcal{U}(X, \textbf{h}, \Lambda, \widetilde{B}) =\bigcap\limits_{i=2}^{n}\ZZ\PP[X_1,X^\prime_1,X^{\pm1}_2,\ldots,X^{\pm1}_{i-1},X_i,X^{\prime}_i, X^{\pm1}_{i+1},\ldots,X^{\pm1}_{n}].
$$
\end{proposition}
\begin{proof}
By the similar argument in \cite[Proposition 4.3]{bfz}, we only need to consider the case when $n=2$, which is proved in Lemma \ref{lemcapn=2}.
\end{proof}

\begin{proposition}\label{propu=u1n=2}
We have
\begin{align*}
&\ZZ\PP[X_1,X^\prime_1,X^{\pm1}_2,\ldots,X^{\pm1}_{j-1},X_j,X^{\prime\prime}_j, X^{\pm1}_{j+1},\ldots,X^{\pm1}_{n}] \\ =&\ZZ\PP[X_1,X^\prime_1,X^{\pm1}_2,\ldots,X^{\pm1}_{j-1},X_j,X^{\prime}_j, X^{\pm1}_{j+1},\ldots,X^{\pm1}_{n}].
\end{align*}
\end{proposition}
\begin{proof}
Without loss of generality, we can assume that $i=1$ and $j=2$.

Case 1. If $b_{0,12}=0$, then $b_{0,21}=b_{1,21}=b_{1,12}=0$. Hence $b^{2}_{1}=b^{2}_{0}$. It follows that
\begin{align*}
X^{\prime\prime}_2=&\sum\limits_{r=0}^{d_2}h_{2,r}(q^{\frac{1}{2}})X^{\prime}(r[\beta_{1}^{2}]_+ +(d_2-r)[-\beta_{1}^{2}]_+ -e_2)\\
=&\sum\limits_{r=0}^{d_2}h_{2,r}(q^{\frac{1}{2}})X(r[\beta_{0}^{2}]_+ +(d_2-r)[-\beta_{0}^{2}]_+ -e_2)=X^{\prime}_2
\end{align*}
and the statement is true.

Case 2. If $b_{0,12}<0$ then $b_{0,21}>0$, $b_{1,12}>0$ and $b_{1,21}<0$. By symmetry, we only need to show that $X^{\prime\prime}_2\in\ZZ\PP[X^{}_{1},X^{\prime}_{1},X^{}_{2},X^{\prime}_{2},X^{\pm1}_3,\ldots,X^{\pm1}_n]$. Recall that
\begin{align}\label{xprimeprime1}
X^{\prime\prime}_2
=\sum\limits_{r=0}^{d_2}h_{2,r}(q^{\frac{1}{2}})X^{\prime}(r\beta_{1}^{2} +[-b_{1}^{2}]_+ -e_2).
\end{align}
Note that $X^{\prime}(r\beta_{1,12}e_1+r\beta_{1}^{2} +[-b_{1}^{2}]_+ -e_2-r\beta_{1,12}e_1) =q^{-\frac{1}{2}\Lambda_1(r\beta_{1,12}e_1,T_r)}(X^{\prime}_1)^{r\beta_{1,12}}X(T_r)$, where $T_r=r\beta_{1}^{2} +[-b^{2}_{1}]_+ -e_2 -r\beta_{1,12}e_1$.

By Lemma \ref{lemxprime}, we have that
\begin{align*}
(X^{\prime}_{1})^{s}
=V_{b_{0}^{1}}^s X(s[-b^{1}_0]_+-se_1),
\end{align*}
where $V_{b_{0}^{1}}^s = \prod\limits_{k=1}^{s} (\sum\limits_{l=0}^{d_1} q^{\frac{2k-1}{2}l\widetilde{d}_1d^{-1}_1} h_{1,l}(q^{\frac{1}{2}})X(l\beta_{0}^{1}) )$ if $s\in\ZZ_+$ and $V_{b_{0}^{1}}^0=1$. It concludes that
\begin{align*}
X^{\prime}(r\beta_{1}^{2} +[-b_{1}^{2}]_+ -e_2)
=q^{-\frac{1}{2}\Lambda_1(r\beta_{1,12} e_1,T_r)}V_{b_{0}^{1}}^{r\beta_{1,12}} X(r\beta_{1,12}[-b^{1}_0]_+-r\beta_{1,12}e_1)X(T_r)
\end{align*}
for $r\in[0,d_2]$. Since $E_{+}^{\widetilde{B}_0}r\beta_{1,12}e_1=r\beta_{1,12}[-b^{1}_0]_+-r\beta_{1,12}e_1$ and $E_{+}^{\widetilde{B}_0}T_r =T_r$, it follows that
$$
\Lambda_1(r\beta_{1,12}e_1,T_r)=\Lambda_0(r\beta_{1,12}[-b^{1}_0]_+-r\beta_{1,12}e_1,T_r)
$$
for $r\in[0,d_2]$. Using the fact that
$
\beta_{1}^{2} -2\beta_{1,12}e_1 +\beta_{1,12}[-b_{0}^{1}]_+ =\beta_{0}^{2},
$
it follows that
\begin{equation}\label{xprimeprime2}
X^{\prime}(r\beta_{1}^{2} +[-b_{1}^{2}]_+ -e_2)=V^{r\beta_{1,12}}_{b_{0}^{1}}X(r\beta_{0}^{2}+[-b_{1}^{2}]_+-e_2)
\end{equation}
for $r\in[0,d_2]$. By (\ref{xprimeprime1}) and (\ref{xprimeprime2}), we conclude that
\begin{equation}\label{xprimeprime3}
X^{\prime\prime}_2=\sum\limits_{r=0}^{d_2}h_{2,r}(q^{\frac{1}{2}})V_{b_{0}^{1}}^{r\beta_{1,12}} X(r\beta^{2}_{0}+[-b^{2}_{1}]_+-e_2).
\end{equation}

It is easy to check that
\begin{equation}\label{equinvertibleelt}
q^{-\frac{1}{2}\widetilde{d}_2}X([b^{2}_{0}]_+) =q^{-\frac{1}{2}\Lambda_0(e_2,[-b^{2}_0]_+)}X_2X^{\prime}_2 -\sum\limits_{r=0}^{d_2-1}q^{-\frac{1}{2}r\widetilde{d}_2d^{-1}_2} h_{2,r}(q^{\frac{1}{2}}) X(r\beta^{2}_{0} +[-b^{2}_{0}]_+). \end{equation}

Note that $q^{-\frac{1}{2}\widetilde{d}_2}X([b^{2}_{0}]_+)$ is an invertible element of $\ZZ\PP[X^{\pm1}_3,\ldots,X^{\pm1}_n]$. To prove that  $X^{\prime\prime}_2\in\ZZ\PP[X^{}_{1},X^{\prime}_{1},X_{2},X^{\prime}_{2},X^{\pm1}_3,\ldots,X^{\pm1}_n]$, we only need to show
$$
X^{\prime\prime}_2q^{-\frac{1}{2}}X([b^{2}_{0}]_+)\in \ZZ\PP[X^{}_{1},X^{\prime}_{1},X^{}_{2},X^{\prime}_{2},X^{\pm1}_3,\ldots,X^{\pm1}_n].
$$

By using (\ref{xprimeprime3}) and (\ref{equinvertibleelt}), we have that
$$
X^{\prime\prime}_2q^{-\frac{1}{2}\widetilde{d}_2}X([b^{2}_{0}]_+)= A_1- A_2 +A_3,
$$
where
\begin{enumerate}
  \item[] $A_1=q^{-\frac{1}{2}\Lambda_0(e_2,[-b^{2}_0]_+)}V_{b_{0}^{1}}^{b_{1,12}} X(b^{2}_{0}+[-b^{2}_{1}]_+-e_2)X_2X^{\prime}_2$,
  \item[] $A_2=\sum\limits_{r=0}^{d_2-1}q^{-\frac{1}{2}r\widetilde{d}_2d^{-1}_2} h_{2,r}(q^{\frac{1}{2}}) V_{b_{0}^{1}}^{b_{1,12}} X(b^{2}_{0}+[-b^{2}_{1}]_+-e_2) X(r\beta^{2}_{0} +[-b^{2}_{0}]_+)$,
  \item[] $A_3=q^{-\frac{1}{2}\widetilde{d}_2}\sum\limits_{r=0}^{d_2-1}h_{2,r}(q^{\frac{1}{2}}) V_{b_{0}^{1}}^{r\beta_{1,12}} X(r\beta^{2}_{0}+[-b^{2}_{1}]_+-e_2)X([b^{2}_{0}]_+)$.
\end{enumerate}
We claim that $A_1,A_2-A_3\in\ZZ\PP[X^{}_{1},X^{\prime}_{1},X^{}_{2},X^{\prime}_{2},X^{\pm1}_3,\ldots,X^{\pm1}_n]$. Indeed, observe that
\begin{align*}
&V_{b_{0}^{1}}^{b_{1,12}} X(b^{2}_{0}+[-b^{2}_{1}]_+-e_2)X_2\\
 =& (X^{\prime}_{1})^{b_{1,12}}X(-b_{1,12}[-b^{1}_0]_+ +b_{1,12}e_1)X(b^{2}_{0}+[-b^{2}_{1}]_+-e_2)X_2.
\end{align*}
We check at once that the first two components of the vector
$$
-b_{1,12}[-b^{1}_0]_+ +b_{1,12}e_1+b^{2}_{0}+[-b^{2}_{1}]_+ -e_2 +e_2
$$
are $(0,0)$, which implies that $X(-b_{1,12}[-b^{1}_0]_+ +b_{1,12}e_1)X(b^{2}_{0}+[-b^{2}_{1}]_+-e_2)X_2$ is a Laurent monomial in $\ZZ\PP[X^{\pm1}_3,\ldots,X^{\pm1}_n]$. It follows that
$$A_1\in\ZZ\PP[X_{1},X^{\prime}_{1},X_{2},X^{\prime}_{2},X^{\pm1}_3,\ldots,X^{\pm1}_n].$$

We now turn to prove that $A_2-A_3\in\ZZ\PP[X^{}_{1},X^{\prime}_{1},X^{}_{2},X^{\prime}_{2},X^{\pm1}_3,\ldots,X^{\pm1}_n]$. Set $t_{1,r}=(d_2-r)\beta_{1,12}$ and $t_{2,r}=r\beta_{1,12}$ for $r\in[0,d_2-1]$. By Lemma \ref{lemxprime}, we have that
\begin{equation}\label{equv1}
V_{b_{0}^{1}}^{t_{2,r}}=(X^{\prime}_1)^{t_{2,r}}X(-t_{2,r}[-b^{1}_{0}]_+ +t_{2,r}e_1),
\end{equation}
\begin{align}\label{equv2}
V_{b_{0}^{1}}^{b_{1,12}}=&(X^{\prime}_1)^{t_{2,r}}(X^{\prime}_1)^{t_{1,r}}X(-t_{1,r}[-b^{1}_{0}]_+ +t_{1,r}e_1) X(-t_{2,r}[-b^{1}_{0}]_+ +t_{2,r}e_1)  \\
=&(X^{\prime}_1)^{t_{2,r}}V_{b_{0}^{1}}^{t_{1,r}}X(-t_{2,r}[-b^{1}_{0}]_+ +t_{2,r}e_1). \nonumber
\end{align}

Note that $X(\beta_{0}^{1})$ is a factor of $V_{b_{0}^{1}}^{t_{1,r}}-1$.

Using (\ref{equ1}), it follows that
\begin{align*}
&X(b^{2}_{0}+[-b^{2}_{1}]_+-e_2)X(r\beta_{0}^{2}+[-b^{2}_{0}]) \\
=&q^{-\frac{1}{2}\widetilde{d}_2 +\frac{1}{2}r\widetilde{d}_2d^{-1}_{2}}X(r\beta^{2}_{0}+[-b^{2}_{1}]_+-e_2)X([b^{2}_0]_+).
\end{align*}
Hence $A_2=q^{-\frac{1}{2}\widetilde{d}_2}\sum\limits_{r=0}^{d_2-1} h_{2,r}(q^{\frac{1}{2}}) V_{b_{0}^{1}}^{b_{1,12}} X(-e_2+r\beta^{2}_{0}+[-b^{2}_{1}]_+)X([b^{2}_0]_+)$.

By using (\ref{equv1}) and (\ref{equv2}), we obatin that
\begin{align*}
&A_2-A_3\\
=&q^{-\frac{1}{2}\widetilde{d}_2}\sum\limits_{r=0}^{d_2-1} h_{2,r}(q^{\frac{1}{2}}) (V_{b_{0}^{1}}^{b_{1,12}} -V_{b_{0}^{1}}^{r\beta_{1,12}}) X(-e_2+r\beta^{2}_{0}+[-b^{2}_{1}]_+)X([b^{2}_0]_+) \\
=& q^{-\frac{1}{2}\widetilde{d}_2}\sum\limits_{r=0}^{d_2-1}h_{2,r}(q^{\frac{1}{2}}) (X^{\prime}_1)^{t_{2,r}} (V^{t_{1,r}}_{b^{1}_0} -1) X(-t_{2,r}[-b^{1}_{0}]_+ +t_{2,r}e_1)\\
&\cdot  X(r\beta^{2}_{0}+[-b^{2}_{1}]_+-e_2)X([b^{2}_{0}]_+).
\end{align*}

For each $r\in[0,d_2-1]$, the first two components of the vector
$$
-t_{2,r}[-b^{1}_{0}]_+ +t_{2,r}e_1 -e_2+r\beta^{2}_{0}+[-b^{2}_{1}]_+ +[b^{2}_{0}]_+
$$
are $(0,-1)$. But the first two components of $\beta_{0}^{1}$ are $(0,\beta_{0,21})$, where $\beta_{0,21}>0$.
Thus we have
$$
(V^{t_{1,r}}_{b^{1}_0} -1) X(-t_{2,r}[-b^{1}_{0}]_+ +t_{2,r}e_1)X(r\beta^{2}_{0}+[-b^{2}_{1}]_+-e_2)X([b^{2}_{0}]_+) \in\ZZ\PP[X^{\pm1}_3,\ldots,X^{\pm1}_n].
$$
Therefore $A_2-A_3\in\ZZ\PP[X^{}_{1},X^{\prime}_{1},X^{}_{2},X^{\prime}_{2},X^{\pm1}_3,\ldots,X^{\pm1}_n]$.

Case 3. If $b_{0,12}>0$ then $b_{0,21}<0$, $b_{1,12}<0$ and $b_{1,21}>0$. Write
$$
X^{\prime\prime}_2
=\sum\limits_{r=0}^{d_2}h_{2,r}(q^{\frac{1}{2}})X^{\prime}(-r\beta_{1}^{2} +[b_{1}^{2}]_+ -e_2)
$$
and $(X^{\prime}_{1})^{s}=W_{b_{0}^{1}}^s X(s[b^{1}_0]_+-se_1)$, where $W_{b^{1}_0}^s = \prod\limits_{k=1}^{s} (\sum\limits_{l=0}^{d_1} q^{\frac{1-2k}{2}l\widetilde{d}_1d^{-1}_1} h_{1,l}(q^{\frac{1}{2}})X(-l\beta^{1}_0) )$ for $s\in\ZZ_+$ and $W_{b^{1}_0}^0=1$.
Note that
$$
X^\prime(-r\beta^{2}_1+[b^{2}_1]_+-e_2)=q^{-\frac{1}{2}\Lambda_1(-r\beta_{1,12}e_1,T^\prime_r)} X^\prime(-r\beta_{1,12}e_1)X(T^\prime_r),
$$
where $T^\prime_r=-r\beta^2_1+[b^2_1]_+-e_2+r\beta_{1,12}e_1$

Because $E^{\widetilde{B}_0}_-T^\prime_r=T^\prime_r$ and $E^{\widetilde{B}_0}_-(-r\beta_{1,12}e_1)=-r\beta_{1,12}([b^1_0]_+-e_1)$, we know that $\Lambda_1(-r\beta_{1,12}e_1,T^\prime_r)=\Lambda_0(-r\beta_{1,12}([b^1_0]_+-e_1),T^\prime_r)$.

By using the fact that $\beta^2_0=\beta^2_1-2\beta_{1,12}e_1+\beta_{1,12}[b^1_0]_+$, we obtain that $X^\prime(-r\beta_{1,12}e_1)=W^{-r\beta_{1,12}}_{b^1_0}X(-r\beta_{1,12}([b^1_0]_+-e_1))$ for $r\in[0,d_2]$. Thus
$$
X^{\prime\prime}_2=\sum\limits_{r=0}^{d_2}h_{2,r}(q^{\frac{1}{2}})W_{b_{0}^{1}}^{-r\beta_{1,12}}X(-r\beta^2_0+[b^2_1]_+-e_2).
$$
By direct calculations, we have that
$$
q^{\frac{1}{2}\widetilde{d}_2}X([-b^2_0]_+)=q^{\frac{1}{2}\Lambda_0([b^2_0]_+,e_2)}X_2X^\prime_2-\sum\limits_{r=0}^{d_{2}-1} q^{\frac{1}{2}r\widetilde{d}_2d^{-1}_2} h_{2,r}(q^{\frac{1}{2}})X(-r\beta^2_0+[b^2_0]_+).
$$
It suffices to show that $X^{\prime\prime}_2q^{\frac{1}{2}\widetilde{d}_2}X([-b^2_0]_+) \in\ZZ\PP[X^{}_{1},X^{\prime}_{1},X^{}_{2},X^{\prime}_{2},X^{\pm1}_3,\ldots,X^{\pm1}_n]$. Similarly,
$$
X^{\prime\prime}_2q^{\frac{1}{2}\widetilde{d}_2}X([-b^2_0]_+)=B_1-B_2+B_3,
$$
where
\begin{align*}
B_1=&q^{\frac{1}{2}\Lambda_0([b^{2}_0]_+,e_2)}W_{b_{0}^{1}}^{-b_{1,12}} X(-b^{2}_{0}+[b^{2}_{1}]_+-e_2)X_2X^{\prime}_2, \\
B_2=&\sum\limits_{r=0}^{d_2-1}q^{\frac{1}{2}r\widetilde{d}_2d^{-1}_2} h_{2,r}(q^{\frac{1}{2}}) W_{b_{0}^{1}}^{-b_{1,12}} X(-b^{2}_{0}+[b^{2}_{1}]_+-e_2) X(-r\beta^{2}_{0} +[b^{2}_{0}]_+) \\
=&q^{-\frac{1}{2}\widetilde{d}_2}\sum\limits_{r=0}^{d_2-1}h_{2,r}(q^{\frac{1}{2}}) W_{b_{0}^{1}}^{-b_{1,12}} X(-r\beta^{2}_{0}+[b^{2}_{1}]_+-e_2)X([-b^{2}_{0}]_+), \\
B_3=&q^{-\frac{1}{2}\widetilde{d}_2}\sum\limits_{r=0}^{d_2-1}h_{2,r}(q^{\frac{1}{2}}) W_{b_{0}^{1}}^{-r\beta_{1,12}} X(-r\beta^{2}_{0}+[b^{2}_{1}]_+-e_2)X([-b^{2}_{0}]_+).
\end{align*}

The rest proof is similar to the one for Case 2, so we omit it. The proof is completed.
\end{proof}

\begin{remark}
The polynomial exchange relations of the generalized quantum cluster algebras make the proof of Proposition \ref{propu=u1n=2} not obviously the same as one in the quantum cluster algebras, which requires more careful discussions
and calculations.

\end{remark}

The following theorem is the main result of this section which can be viewed as the quantum analogue of \cite[Theorem 4.1]{gsv18}.
\begin{theorem}\label{thmu=u1}
If the quantum seed $(X,\mathbf{h},\Lambda_0,\widetilde{B}_0)$ is coprime, then
$$
\mathcal{U}(X, \textbf{h}, \Lambda_0, \widetilde{B}_0)= \mathcal{U}(X^{\prime},\mathbf{h},\Lambda_1,\widetilde{B}_1).
$$
Moreover, if every quantum seed, which is mutation-equivalent to the initial seed $(X,\mathbf{h},\Lambda_0,\widetilde{B}_0)$, is coprime, then any upper bound of $\A(X,\mathbf{h},\Lambda_0,\widetilde{B}_0)$ coincides with the generalized quantum upper cluster algebra $\widetilde{\mathcal{U}}(X,\mathbf{h},\Lambda_0,\widetilde{B}_0)$.
\end{theorem}
\begin{proof}
Assume that $i\neq j$ and let $X^{\prime\prime}_j=\mu_j\mu_i(X_j)$. By Proposition \ref{propu2}, we have that $\mathcal{U}(X^\prime,\textbf{h},\Lambda_1,\widetilde{B}_1) =\bigcap\limits_{j=2}^{n}\ZZ\PP[X_1,X^\prime_1,X^{\pm1}_2,\ldots,X^{\pm1}_{j-1},X_j,X^{\prime\prime}_j, X^{\pm1}_{j+1},\ldots,X^{\pm1}_{n}]$ and by Proposition \ref{propu=u1n=2}, it follows that $\mathcal{U}(X,\textbf{h},\Lambda_0,\widetilde{B}_0)=\mathcal{U}(X^\prime,\textbf{h},\Lambda_1,\widetilde{B}_1)$.
\end{proof}

\begin{remark}
Theorem \ref{thmu=u1} implies that under the ``coprimality" condition, the  generalized quantum cluster algebras of geometric types have the Laurent phenomenon.
\end{remark}



\begin{thebibliography}{99}

\bibitem{BCDX}
L. Bai, X. Chen, M. Ding and F. Xu.
\emph{A quantum analogue of generalized cluster algebras}. Algebr. Represent. Theory 21 (2018), no.~6, 1203--1217.

\bibitem{BCDX-1}
L. Bai, X. Chen, M. Ding and F. Xu.
\emph{On the generalized cluster algebras of geometric types.} Symmetry Integrability Geom. Methods Appl. 16 (2020), Paper No. 092.

\bibitem{bfz}
A. Berenstein, S. Fomin and A. Zelevinsky.
\emph{Cluster algebras III: Upper bounds and double Bruhat cells.} Duke Math. J. \textbf{126} (2005), no.~1, 1--52.

\bibitem{BZ2005}
A. Berenstein and A. Zelevinsky. \emph{Quantum cluster algebras}. Adv. Math. \textbf{195} (2005), no. 2, 405--455.






\bibitem{CL1} P. Cao and F.  Li. \emph{Some conjectures on generalized cluster algebras via the cluster formula and D-matrix pattern}. J. Algebra 493 (2018), 57--78.

\bibitem{CL2} P. Cao and F. Li. \emph{On some combinatorial properties of generalized cluster algebras}. J. Pure Appl. Algebra 225 (2021), no. 8, Paper No. 106650, 13 pp.




\bibitem{CS}
L.~Chekhov and M.~Shapiro. \emph{Teichm\"uller spaces of Riemann surfaces with orbifold
points of arbitrary order and cluster variables.} Int. Math. Res. Notices 2014, no.~10,
2746--2772.


\bibitem{DL} Q.~Du and F.~Li. \emph{Some elementary properties of Laurent phenomenon algebras.} arXiv:2201.02917.


\bibitem{ca1}
S.~Fomin and A.~Zelevinsky. \emph{Cluster algebras. I. Foundations.}
J. Amer. Math. Soc.  \textbf{15}  (2002),  no. 2, 497--529.

\bibitem{ca2}
S.~Fomin and A.~Zelevinsky. \emph{Cluster algebras. II. Finite type
classification.} Invent. Math.  \textbf{154}  (2003),  no. 1,
63--121.


\bibitem{gsv18}
M. Gekhtman, M. Shapiro and A. Vainshtein. \emph{Drinfeld double of $GL_n$ and generalized cluster structure.}
Proc. Lond. Math. Soc. (3) \textbf{116} (2018), no. ~3, 429--484.



\bibitem{gy}
K. Goodearl and M. Yakimov. \emph{Quantum cluster algebra structures on quantum nilpotent algebras.} Mem. Amer. Math. Soc. \textbf{247} (2017), no.~1169, vii+119 pp.





\bibitem{Mou} L. Mou. \emph{Scattering diagrams for generalized cluster algebras}. arXiv:2110.02416.






\bibitem{nak2}
T.~Nakanishi. \emph{Structure of seeds in generalized cluster algebras.} Pacific J. Math. \textbf{277} (2015), no.~1,
201--217.

\end{thebibliography}
\end{document}